\RequirePackage{fix-cm}

\documentclass[leqno, fleqn, centertags, 12pt]{article}

\usepackage{latexsym}
\usepackage{amsmath}
\usepackage{amssymb}
\usepackage{verbatim}

\usepackage[T1]{fontenc}

\usepackage[upright]{fourier}

\usepackage{bm}

\usepackage{fullpage}

\usepackage{tikz-cd}

\tikzcdset{row sep/special/.initial=12ex}
\tikzcdset{row sep/spe/.initial=10.5ex}
\tikzcdset{row sep/speh/.initial=13ex}

\tikzcdset{column sep/spec/.initial=5em}
\tikzcdset{column sep/spech/.initial=6em}
\tikzcdset{column sep/speci/.initial=6em}
\tikzcdset{column sep/specc/.initial=8em}
\tikzcdset{column sep/speccc/.initial=10em}

\tikzcdset{row sep/red/.initial=3.79em}

\tikzcdset{row sep/tri/.initial=5em}
\tikzcdset{column sep/tri/.initial=2.4em}

\tikzcdset{row sep/boom/.initial=5em}
\tikzcdset{column sep/boom/.initial=5.1em}

\tikzcdset{row sep/boomm/.initial=5.6em}
\tikzcdset{column sep/boomm/.initial=4.25em}

\tikzcdset{row sep/boommm/.initial=6.32em}
\tikzcdset{column sep/boommm/.initial=6.45em}

\tikzcdset{row sep/boommmm/.initial=7em}
\tikzcdset{column sep/boommmm/.initial=7.14em}

\tikzcdset{row sep/bboom/.initial=8em}
\tikzcdset{column sep/bboom/.initial=8em}

\tikzcdset{row sep/booms/.initial=4.25em}
\tikzcdset{column sep/booms/.initial=4.335em}

\tikzcdset{row sep/boomss/.initial=4.75em}
\tikzcdset{column sep/boomss/.initial=4.845em}

\tikzcdset{row sep/free/.initial=12ex}
\tikzcdset{column sep/freee/.initial=4em}

\usepackage{comment}

\newskip\nineskipamount \nineskipamount=9pt plus 0pt minus 0pt
\newskip\zeroskipamount \zeroskipamount=0pt plus 0pt minus 0pt
\usepackage[noorphans,vskip=\nineskipamount]{quoting}

\makeatletter
\renewcommand{\@makefntext}[1]{\vspace*{0.5ex}\parindent=0em
\hspace*{-0.4em}
\hbox to 0.4em{\hss\@makefnmark}\hspace*{0.4em}{#1}
}
\makeatother

\newcounter{mysectionnumber}
\newcommand{\mysection}[2]{\setcounter{footnote}{0}
\setcounter{equation}{0}
\setcounter{myparnum}{0}
\refstepcounter{mysectionnumber}
\vspace{27pt}{\Large {\themysectionnumber.} {#1}}\label{#2}\vspace*{15pt}}

\numberwithin{equation}{section}

\newcommand{\myuppar}[1]{\vspace{\medskipamount}\textbf{#1}\hspace*{0.5em}}

\newcommand{\myit}[1]{\textbf{\textit{#1}}\hspace{0.0em}}

\newcounter{myparnum}[mysectionnumber]
\renewcommand{\themyparnum}{\arabic{mysectionnumber}.\arabic{myparnum}}
\newcommand{\mypar}[2]{\refstepcounter{myparnum}{\vspace{\medskipamount}\textbf{{\themyparnum. #1}\label{#2}}\hspace{0.5em}}}

\newcounter{mylemmanum}[myparnum]
\newcounter{myaparnum}
\newcommand{\myappend}[2]{\setcounter{footnote}{0}
\setcounter{myaparnum}{0}
\vspace{27pt}{\Large A\dff.\oss {#1}}\label{#2}\vspace*{15pt}}

\newcommand{\myapar}[2]{\refstepcounter{myaparnum}{\vspace{\medskipamount}\textbf{{\themyaparnum. #1}\label{#2}}\hspace{0.5em}}}
\renewcommand{\themyaparnum}{A\halfff\fff.\fff\arabic{myaparnum}}

\newcommand{\proof}{\vspace{\medskipamount}{\textbf{{\emph{Proof}.}}\hspace*{1em}}}

\newcommand{\eproof}{ $\blacksquare$}

\newcommand{\dis}{\displaystyle}

\def\sss{\hspace{0.05em}\ }
\def\dss{\hspace{0.1em}\ }
\def\trs{\hspace{0.15em}\ }
\def\qss{\hspace{0.2em}\ }
\def\pss{\hspace{0.3em}\ }
\def\oss{\hspace{0.4em}\ }

\def\halfff{\hspace*{0.025em}}
\def\fff{\hspace*{0.05em}}
\def\dff{\hspace*{0.1em}}
\def\trf{\hspace*{0.15em}}
\def\qff{\hspace*{0.2em}}
\def\pff{\hspace*{0.3em}}
\def\off{\hspace*{0.4em}}

\def\ttff{{\hspace*{-0.05em}--\hspace*{0.15em}}}

\newcommand{\nsp}{\hspace*{-0.1em}}
\newcommand{\nnsp}{\hspace*{-0.15em}}

\newcommand{\dnsp}{\hspace*{-0.2em}}

\renewcommand{\leq}{\leqslant}
\renewcommand{\geq}{\geqslant}

\newcommand{\supp}{\mathop{\mbox{\textup{supp}}}}

\newcommand{\rrr}{\mathbf{R}}
\newcommand{\nnn}{\mathbf{N}}

\newcommand{\inte}{\operatorname{int}\trf}

\newcommand{\sub}{\operatorname{\textit{sub}}\qff}

\newcommand{\cf}{{\dff\operatorname{c{\fff}f}}}
\newcommand{\lf}{{\dff\operatorname{l{\halfff}f}}}
\newcommand{\iinf}{{\dff\operatorname{inf}}}
\newcommand{\lone}{{\dff l_{\dff 1}}}

\newcommand{\ry}{r_{\dff \smallsetminus\dff Y}}

\newcommand{\hry}{r_{\dff \smallsetminus\dff Y\trf *}}

\newcommand{\num}[1]{|\qff #1 \qff|}

\newcommand{\norm}[1]{\|\qff #1 \qff\|}

\newcommand{\ssm}[0]{\qff \smallsetminus\qff}

\newcommand{\ttoo}{\hspace*{0.2em}\longrightarrow\hspace*{0.2em}}

\newcommand{\particular}{par\-tic\-u\-lar}

\begin{document}

\setlength{\baselineskip}{12pt plus 0pt minus 0pt}
\setlength{\parskip}{12pt plus 0pt minus 0pt}
\setlength{\abovedisplayskip}{12pt plus 0pt minus 0pt}
\setlength{\belowdisplayskip}{12pt plus 0pt minus 0pt}

\newskip\smallskipamount \smallskipamount=3pt plus 0pt minus 0pt
\newskip\medskipamount   \medskipamount  =6pt plus 0pt minus 0pt
\newskip\bigskipamount   \bigskipamount =12pt plus 0pt minus 0pt

\author{Nikolai\qss V.\qss Ivanov}
\title{Leray\qss theorems\qss for\dss $l_{\fff 1}$\dnsp-norms\qss of\pss infinite\qss chains}
\date{}

\footnotetext{\hspace*{-0.65em}\copyright\oss 
Nikolai\qss V.\qss Ivanov,\oss 2020.\trs 
Neither\sss the work reported\sss in\dss the present paper\halfff,\qss
nor\dss its preparation were supported\dss by\dss any\sss corporate entity.}

\maketitle

\renewcommand{\baselinestretch}{1}
\selectfont

\vspace*{12ex}

\myit{\hspace*{0em}\large Contents}\vspace*{1ex} \vspace*{\bigskipamount}\\ 
\hbox to 0.8\textwidth{\myit{\phantom{A.}1.}\hspace*{0.5em} Introduction\hfil  2}\hspace*{0.5em} \vspace*{0.25ex}\\
\hbox to 0.8\textwidth{\myit{\phantom{A.}2.}\hspace*{0.5em} A\dss Leray\trs theorem\dss for\dss infinite chains\hfil 6}\hspace*{0.5em} \vspace*{0.25ex}\\
\hbox to 0.8\textwidth{\myit{\phantom{A.}3.}\hspace*{0.5em} Compactly\dss finite and $l_{\dff 1}$\dnsp-homology\hfil 12}\hspace*{0.5em} \vspace*{0.25ex}\\
\hbox to 0.8\textwidth{\myit{\phantom{A.}4.}\hspace*{0.5em} Extensions of\dss coverings and $l_{\dff 1}$\dnsp-homology\hfil 19}\hspace*{0.5em} \vspace*{0.25ex}\\
\hbox to 0.8\textwidth{\myit{\phantom{A.}5.}\hspace*{0.5em} Removing\dss weakly $l_{\dff 1}$\dnsp-acyclic subspaces\hfil 20}\hspace*{0.5em} \vspace*{1ex}\\
\hbox to 0.8\textwidth{\myit{Appendix.}\hspace*{0.5em} Double complexes\hfil 26}\hspace*{0.5em}
\vspace*{1ex}\\
\hbox to 0.8\textwidth{\myit{References}\hspace*{0.5em}\hfil 27}\hspace*{0.5em}  \vspace*{0.25ex}

\renewcommand{\baselinestretch}{1}
\selectfont

\vspace*{21ex}

{\small
The paper\dss is\dss devoted\dss to an adaptation of\dss author's approach\qss \cite{i3}\qss
to\trs Leray\trs theorems\sss in\sss bounded cohomology\dss theory\dss to\sss infinite chains.\oss 
The paper\sss may\sss be considered as a continuation of\trs the paper\qss \cite{i3},\oss
but\sss depends on\sss it\dss mostly\dss for\sss the motivation of\dss proofs.\oss
Among\sss the main results are a stronger and\sss more general\dss form of\trs Gromov's\trs 
Vanishing-finiteness\sss theorem and a generalization of\trs the first\sss part\sss 
of\trs his\dss Cutting-of\trs theorem.\oss 
The proofs do not\sss depend on any\sss tools specific for\sss
the bounded cohomology and $l_{\dff 1}$\dnsp-homology\sss theory,\oss 
but\sss use\sss the fact\dss that  $l_{\dff 1}$\dnsp-homology depend only\sss
on\sss the fundamental\dss group.\oss 
}

\newpage
\mysection{Introduction}{introduction}

\myuppar{Locally,\pss compactly,\pss and\sss star\sss finite families.}
A\qss \emph{family\sss of\dss subsets}\qss
$\mathcal{U}
\off =\off
\{\qff
U_{\dff i}
\qff\}_{\dff i\dff \in\qff I}$\dss
of\dss a set\sss $X$\sss is\dss a map\sss
$i\off \longmapsto\off U_{\dff i}\pff \subset\off X$\sss
from a set\sss $I$\sss to\sss the set\sss of\dss all\sss subsets of\dss $X$\nnsp.\oss
The\sss family\sss 
$\mathcal{U}$\dss
is\dss said\dss to be\qss 
\emph{star\sss finite}\pss if\dss for every $i\qff \in\pff I$\sss
the intersection\sss 
$U_{\dff i}\qff \cap\qff U_{\dff j}$\sss is\dss non-empty\sss
for only\sss a\sss finite number of\dss $j\qff \in\pff I$\nnsp.\oss
Usually,\oss but\dss not\sss always,\qss only\sss the set
$\{\qff U_{\dff i}\qff |\qff i\qff \in\qff I\pff\}$\sss matters,\pss
and\sss we write $U\qff \in\qff \mathcal{U}$ instead of\pss
``\nnsp$U\off =\off U_{\dff i}$\sss for some\sss $i\qff \in\pff I$\nnsp''.\oss

Let $X$ be a\sss topological\sss space and
$\mathcal{U}
\off =\off
\{\qff
U_{\dff i}
\qff\}_{\dff i\dff \in\qff I}$
be a family\sss of\dss subsets of\dss $X$\nnsp.\oss
The\sss family\sss 
$\mathcal{U}$
is\dss said\dss to be\qss 
\emph{locally\dss finite}\pss if\dss for every $x\qff \in\qff X$
there exists an open neighborhood $V$ of\dss $x$ such\dss that $x\qff \in\qff U$
and $V\qff \cap\qff U_{\dff i}\off \neq\off \varnothing$
for only\sss a\sss finite number of\dss $i\qff \in\pff I$\nnsp,\oss
and\qss \emph{compactly\dss finite}\pss if\dss for every\sss compact\sss 
$Z\qff \subset\qff X$\sss the intersection\sss 
$Z\qff \cap\qff U_{\dff i}\off \neq\off \varnothing$\sss
for only\sss a\sss finite number of\dss $i\qff \in\pff I$\nnsp.\oss
For\sss locally\sss compact\sss spaces\sss the notions of\trs
locally\sss finite and compactly\sss finite families are equivalent.\oss
Eventually\sss our assumptions will\dss imply\sss that\sss $X$\sss is\dss locally\sss compact,\oss
but\sss we prefer\sss to be precise about\sss which\sss finiteness condition\dss is\dss used.\oss
Gromov\qss \cite{gro},\pss  L\"{o}h\dss and\dss Sauer\qss \cite{ls},\pss 
and\dss Frigerio\dss and\dss Moraschini\qss \cite{fm}\qss 
call\sss compactly\dss finite\sss families\qss ``locally\dss finite''.\oss

\myuppar{Coverings.}
The family\sss $\mathcal{U}$\sss is\dss said\dss to be a\qss \emph{covering}\qss of\dss $X$\sss
if\trs the union\dss $\cup_{\dff i\dff \in\qff I}\qff U_{\dff i}$\dss is\dss
equal\dss to\sss $X$\nnsp.\oss
For a covering\sss $\mathcal{U}$\sss we will\sss denote by\sss 
$\mathcal{U}^{\dff \cap}$\sss the collection of\dss all\sss
non-empty\sss finite intersection of\dss elements of\dss $\mathcal{U}$\dnsp.\oss
A covering\sss $\mathcal{U}$\sss is\dss said\dss to be\qss \emph{open}\pff
if\dss every\sss $U_{\dff i}$\sss is\dss open,\oss and\qss 
\emph{proper}\pss if\trs the interiors\sss $\inte U_{\dff i}$\sss
form a covering\sss of\dss $X$\sss and\dss the closures of\trs the sets\sss $U_{\dff i}$\sss
are compact.\oss
Clearly,\oss if\qss there exists a proper covering\sss of\dss $X$\nnsp,\oss
then\dss $X$\dss is\dss locally\sss compact.\oss
It\dss is\dss easy\dss to see\sss that\sss a proper covering\dss
is\dss compactly\sss finite\dss if\trs and\dss only\dss if\trs
it\dss is\dss star finite.\oss

\myuppar{Locally\sss and compactly\sss finite singular\sss chains and\dss homology.}
Recall\dss that\sss a singular $n$\dnsp-sim\-plex\sss in\sss $X$\sss
is\dss a continuous map\sss
$\sigma\dff \colon\dff
\Delta^n\qff \ttoo\qff X$\nnsp.\oss
Let\dss $S_{\fff n}\dff(\trf X\trf)$\sss be\sss the set\sss of\dss singular
$n$\dnsp-simplices in\dss $X$\nnsp.\oss
A subset\dss $I\qff \subset\qff S_{\fff n}\dff(\trf X\trf)$\sss
is\dss said\dss to be\qss \emph{locally\dss finite}\pss
if\trs the family\sss of\dss images\dss 
$\{\qff \sigma\dff(\dff \Delta^n\trf)\qff\}_{\qff \sigma\dff \in\qff I}$\sss
is\trs locally\sss finite,\oss
and\qss \emph{compactly\dss finite}\pss if\trs this family\dss is\dss
compactly\sss finite.\oss
An\qss \emph{infinite\dss singular $n$\dnsp-chain}\pss in\sss $X$\sss 
is\dss defined as a\sss formal\sss sum\vspace{3pt}
\begin{equation}
\label{singular-chain}
\quad
c
\off =\off
\sum\nolimits_{\qff \sigma\qff \in\qff S_{\fff n}\dff(\trf X\trf)}\qff
a_{\dff \sigma}\dff \sigma
\end{equation}

\vspace{-9pt}
with coefficients\sss $a_{\dff \sigma}\qff \in\qff A$\nnsp,\oss
where\sss $A$\sss is\dss some abelian\sss group.\oss
Let $C_{\dff n}^\iinf\dff(\trf X\fff,\qff A\trf)$\sss
be\sss the group of\dss such chains.\oss
The chain $c$\sss
is\dss said\dss to be\qss \emph{locally\dss finite}\pss
if\trs\vspace{3pt}
\[
\quad
\mathcal{S}_{\dff c}
\off =\off
\bigl\{\pff \sigma\dff(\trf \Delta^n\trf)
\pff \bigl|\pff
\sigma\qff \in\qff S_{\fff n}\dff(\trf X\trf)\fff,\pff
a_{\dff \sigma}\off \neq\off 0 \pff\bigr\}
\]

\vspace{-9pt}
is\trs locally\dss finite,\oss
and\qss \emph{compactly\dss finite}\pss if\dss $\mathcal{S}_{\dff c}$\sss is\dss
compactly\sss finite.\oss
The groups of\dss locally\sss and\sss compactly\sss
finite chains are denoted\dss by\sss
$C_{\dff n}^\lf\dff(\trf X\fff,\qff A\trf)$\sss
and\sss
$C_{\dff n}^\cf\dff(\trf X\fff,\qff A\trf)$\sss
respectively.\oss
It\dss is\dss easy\dss to see\sss that\sss
every\dss locally\dss finite chain\dss is\dss compactly\dss finite.\oss
If\dss $c$\sss is\dss compactly\sss finite chain,\oss
then\sss in\dss the usual\sss formula for\dss the boundary\dss
$\partial\dff c$\sss the coefficient\sss in\sss front\sss of\dss each singular
simplex\dss is\dss a\sss finite sum.\oss
Therefore\sss the boundaries\sss $\partial\dff c$\sss of\dss 
compactly\dss finite chains,\oss and\dss hence of\trs locally\dss finite chains,\oss 
are\sss well\sss defined.\oss
Since a singular simplex\sss has only\sss finite number of\dss faces,\oss
the boundary\sss of\dss a\sss locally\sss finite
singular chain\dss is\dss locally\sss finite,\oss
and\dss the boundary\sss of\dss a\sss compactly\sss finite singular chain\dss
is\dss compactly\sss finite.\oss
This\sss leads\sss to\sss two\sss types of\dss singular\sss homology\dss groups
based on\sss in\sss infinite chains,\oss
namely\dss the homology\dss groups\sss
$H_{\dff *}^\lf\dff(\trf X\fff,\qff A\trf)$\sss
based on\sss locally\dss finite chains,\oss
and\dss the homology\dss groups\sss
$H_{\dff *}^\cf\dff(\trf X\fff,\qff A\trf)$\sss
based on\sss compactly\dss finite chains.\oss
From now on we will\sss assume\sss that\sss
$A\off =\off \rrr$\sss and\sss omit\dss the coefficient\sss group.\oss

\myuppar{The norms of\trs infinite singular\sss chains.}
The $l_{\dff 1}$\dnsp-norm\sss $\norm{c}$\sss of\trs the singular chain\qss
(\ref{singular-chain})\qss is\dss\vspace{3pt}
\[
\quad
\norm{c}
\off =\off
\sum\nolimits_{\qff \sigma\qff \in\qff S_{\fff n}\dff(\trf X\trf)}\qff
\num{a_{\dff \sigma}}
\pff.
\]

\vspace{-9pt}
It\dss may\sss happen\dss that\dss
$\norm{c}
\off =\off
\infty$\nnsp.\oss
The $l_{\dff 1}$\dnsp-norm\sss $\norm{h}$ of\dss a homology\sss class\dss
$h\qff \in\pff H_{\dff *}^\cf\trf(\trf X\trf)$\sss
or\sss 
$h\qff \in\pff H_{\dff *}^\lf\trf(\trf X\trf)$\sss 
is\dss defined as\dss 
$\norm{h}
\off =\off\dff 
\inf\pff \norm{c}$\nnsp,\oss
where\sss the infimum\dss is\dss taken over\sss all\sss chains $c$
representing\sss the homology class $h$\nnsp.\oss
Again,\oss it\dss may\sss happen\dss that\dss
$\norm{h}
\off =\off
\infty$\nnsp.\oss

\myuppar{Singular $l_{\dff 1}$\dnsp-homology.}
For\sss an\sss integer\sss $n\qff \geq\qff 0$\dss let\sss
$L_{\trf n}\fff(\trf X\trf)$\sss be\sss the vector space\sss
of\dss infinite singular $n$\dnsp-chains with\sss real\sss coefficients
having\dss finite $l_{\dff 1}$\dnsp-norm.\oss
Such chains are called\dss \emph{$l_{\dff 1}$\dnsp-chains}\pss of\dss dimension $n$\nnsp.\oss
The $l_{\dff 1}$\dnsp-norm\dss turns\sss $L_{\trf n}\fff(\trf X\trf)$\sss
into a\dss Banach\dss space.\oss
The vector space\dss $C_{\trf n}\fff(\trf X\trf)$\dss 
of\dss finite singular $n$\dnsp-chains in\sss $X$\sss
is\dss dense in\dss $L_{\trf n}\fff(\trf X\trf)$\dss
and\dss the boundary\sss operator\sss\vspace{3pt}
\[
\quad
\partial\dff \colon\dff
C_{\trf n}\fff(\trf X\trf)
\qff \ttoo\qff
C_{\trf n\dff -\dff 1}\fff(\trf X\trf)
\]

\vspace{-9pt}
extends by continuity\dss to\sss a map\qss
$\partial\dff \colon\dff
L_{\trf n}\fff(\trf X\trf)
\qff \ttoo\qff
L_{\trf n\dff -\dff 1}\fff(\trf X\trf)$\nnsp,\oss
also called\dss the\qss \emph{boundary\dss operator}.\oss
These boundary\sss operators\sss turn\sss $L_{\trf \bullet}\fff(\trf X\trf)$\sss
into a chain complex.\oss
The homology of\trs is\dss complex are knows as\dss
\emph{$l_{\dff 1}$\dnsp-homology}\qss of\dss $X$\sss 
and are denoted\dss by\sss $H_{\dff *}^{\lone}\fff(\trf X\trf)$\nnsp.\oss
The real\sss vector spaces\sss $H_{\dff n}^{\lone}\fff(\trf X\trf)$\sss
inherit $l_{\dff 1}$\dnsp-norms from\dss $L_{\trf n}\fff(\trf X\trf)$\nnsp,\oss
but\sss in\sss general\sss are not\trs Banach\dss spaces,\oss
because non-zero $l_{\dff 1}$\dnsp-homology\sss classes 
may\sss have $l_{\dff 1}$\dnsp-norm equal\dss to $0$\nnsp.\oss

\myuppar{Acyclicity of\dss subsets.}
As\sss in\qss \cite{i3},\oss let\dss us\sss call\sss a\sss topological\sss space\sss $X$\qss
\emph{boundedly\dss acyclic}\pss if\dss its bounded cohomology\sss are isomorphic\sss
to\sss the bounded cohomology\sss of\dss a point.\oss
This property\dss is\dss equivalent\dss to $X$ being\sss path connected
and\sss its fundamental\dss group being boundedly\sss acyclic in an obvious sense.\oss
By\sss a\sss theorem of\trs Sh.\dss Matsumoto\dss and\qss Sh.\dss Morita\qss \cite{mm}\qss
the space\sss $X$\sss is\dss boundedly\sss acyclic\dss if\trs and\dss only\trs if\trs
it\dss is\dss path connected and\dss 
$H_{\dff *}^{\lone}\fff(\trf X\trf)\off =\off 0$\sss
for\sss $n\qff \geq\qff 1$\nnsp.\oss
See also\qss \cite{i3},\oss Theorem\qss 5.1\qss for a proof.\oss
In\dss this paper\sss we are dealing only\sss with\sss homology\sss
and\sss will\sss call\dss boundedly\sss acyclic spaces and\dss groups\dss
\emph{$l_{\dff 1}$\dnsp-acyclic}.\oss

In\qss \cite{i3}\qss a path connected\sss subset\sss $Z$ of\dss $X$\sss
was called\qss \emph{weakly\dss boundedly\sss acyclic}\pss
if\trs the image of\trs the inclusion\sss homomorphism\dss
$\pi_{\dff 1}\dff(\trf Z\fff,\qff z \trf)
\trf \ttoo\trf 
\pi_{\dff 1}\dff(\trf X\fff,\qff z \trf)$\dss
is\dss boundedly\sss acyclic,\oss
i.e.\qss is\sss $l_{\dff 1}$\dnsp-acyclic in our current\dss terminology.\oss
The ambient\sss space $X$ was fixed.\oss
Now\sss we need a more flexible\sss version of\trs this notion.\oss
Suppose\sss that\sss $Z\qff \subset\qff Y\qff \subset\qff X$\nnsp.\oss
The subset\sss $Z$\sss is\dss said\dss to be\qss
\emph{weakly $l_{\dff 1}$\dnsp-acyclic\sss in}\sss $Y$\trs
if\trs the image of\trs the homomorphism\dss
$\pi_{\dff 1}\dff(\trf Z\fff,\qff z \trf)
\trf \ttoo\trf 
\pi_{\dff 1}\dff(\trf Y\fff,\qff z \trf)$\trs
is\sss $l_{\dff 1}$\dnsp-acyclic.\oss

\myuppar{Acyclicity of\dss families and coverings.}
Let $\mathcal{U}$ be a family\sss of\dss subsets\sss of\dss $X$\nnsp.\oss
It\dss is\dss said\dss to be\dss \emph{$l_{\dff 1}$\dnsp-acyclic}\pss 
if\dss every\sss $U\qff \in\qff \mathcal{U}$\sss is\sss $l_{\dff 1}$\dnsp-acyclic.\oss
A covering $\mathcal{U}$ of\dss $X$\sss is\dss said\dss to be\sss 
\emph{$l_{\dff 1}$\dnsp-acyclic}\qss (as a covering\fff)\qss if\trs the family\sss
$\mathcal{U}^{\dff \cap}$\sss is\sss $l_{\dff 1}$\dnsp-acyclic.\oss

A\sss family\sss $\mathcal{U}$\sss is\dss said\sss to be\qss
\emph{almost\sss $l_{\dff 1}$\dnsp-acyclic}\pss 
if\dss every\sss
$U\qff \in\qff \mathcal{U}$\sss is\sss $l_{\dff 1}$\dnsp-acyclic,\oss
except\halfff,\pss perhaps,\oss of\dss a single exceptional\sss element\sss 
$U_{\fff e}\qff \in\qff \mathcal{U}$\dnsp.\oss
A covering $\mathcal{U}$ of\dss $X$\sss is\dss said\dss to be\qss 
\emph{almost\sss $l_{\dff 1}$\dnsp-acyclic}\qss if\trs the family\sss
$\mathcal{U}^{\dff \cap}$\sss is\dss almost\sss $l_{\dff 1}$\dnsp-acyclic
and\sss the exceptional\sss set\sss $U_{\fff e}$ belongs\sss to\sss $\mathcal{U}$\dnsp.\oss

We need an analogue of\dss weakly\sss boundedly\sss acyclic coverings from\qss \cite{i3}.\oss
Requiring\sss sets\sss $U\qff \in\qff \mathcal{U}$\sss
to be weakly $l_{\dff 1}$\dnsp-acyclic in\sss $X$\sss is\dss not\sss sufficient\dss
for working\sss with compactly\sss finite chains.\oss

A family $\mathcal{U}$ of\dss subsets\sss of\dss $X$\sss 
is\dss said\dss to be\qss \emph{weakly $l_{\dff 1}$\dnsp-acyclic}\pss 
if\dss for every\sss $U\qff \in\qff \mathcal{U}$\sss
a subset\sss $U_{\dff +}\qff \subset\qff X$\sss is\dss given,\oss
such\dss that\sss $U\qff \subset\qff U_{\dff +}$\nsp,\oss
the set\sss $U$\sss is\dss weakly $l_{\dff 1}$\dnsp-acyclic\sss in\sss $U_{\dff +}$\nsp,\oss
and\dss the family\sss of\dss subsets\sss $U_{\dff +}$\sss
is\dss compactly\sss finite.\oss
A\sss finite number of\dss subsets\sss $U_{\dff +}$\sss 
can\sss be equal\dss to\sss $X$\nnsp.\oss
A covering $\mathcal{U}$ of\dss $X$\sss is\dss said\dss to be\qss
\emph{weakly $l_{\dff 1}$\dnsp-acyclic}\pss if\trs the family\sss
$\mathcal{U}^{\dff \cap}$\sss is\dss weakly $l_{\dff 1}$\dnsp-acyclic.\oss

Similarly,\oss a\sss family\sss $\mathcal{U}$\sss is\dss said\sss to be\qss
\emph{almost\dss weakly $l_{\dff 1}$\dnsp-acyclic}\pss 
if\dss subsets\sss $U_{\dff +}\qff \subset\qff X$\sss
with\dss the properties\sss listed\sss above 
are given\sss for\sss every subset\sss
$U\qff \in\qff \mathcal{U}$\dnsp,\oss
except\halfff,\pss perhaps,\oss of\dss a single exceptional\sss set\sss 
$U_{\fff e}\qff \in\qff \mathcal{U}$\dnsp.\oss
A covering $\mathcal{U}$ of\dss $X$\sss is\dss said\dss to be\qss 
\emph{almost\sss weakly\sss $l_{\dff 1}$\dnsp-acyclic}\qss if\trs the family\sss
$\mathcal{U}^{\dff \cap}$\sss is\dss almost\sss weakly\sss $l_{\dff 1}$\dnsp-acyclic
and\sss the exceptional\sss set\sss $U_{\fff e}$ belongs\sss to\sss $\mathcal{U}$\dnsp.\oss

\myuppar{Infinite chains\sss in\sss simplicial\sss complexes.}
A simplicial\sss complex $S$\sss is\dss said\dss to be\qss
\emph{star\sss finite}\pss if\dss each\sss its simplex\dss is\dss contained\sss
in\sss only\sss a finite number of\dss simplices.\oss
Equivalently,\pss $S$\sss is\qss
\emph{star\sss finite}\pss 
if\trs the family\sss of\dss its\sss simplices\dss is\dss star\sss finite.\oss
If\dss $\mathcal{U}$\sss is\dss a  family\sss of\dss subsets of\dss $X$\nnsp,\oss
then\sss $\mathcal{U}$\sss is\dss star\sss finite\dss if\trs and\dss only\trs if\trs
the nerve\sss $N_{\trf \mathcal{U}}$\sss of\dss $\mathcal{U}$\sss is\dss star\sss finite.\oss

An\qss \emph{infinite $n$\dnsp-chain}\qss of\dss a simplicial\sss complex $S$\sss
is\dss a potentially\sss infinte formal\sss sum of\sss $n$\dnsp-sim\-plices of\dss $S$\sss
with coefficients in some abelian\sss group\sss $A$\nnsp.\oss
If\dss $S$\sss is\dss star\sss finite,\oss then\sss the usual\sss formula
defines\sss the boundaries\dss $\partial\dff c$\dss 
of\dss infinite $n$\dnsp-chains of\dss $S$\nnsp.\oss
This\sss leads\sss to homology\dss groups\sss
$H_{\dff *}^\iinf\dff(\trf S\fff,\qff A\trf)$\sss
based on\sss infinite chains.\oss
We will\sss always assume\sss that\sss $A\off =\off \rrr$\nnsp.\oss\vspace{-0.125pt}

\myuppar{Leray\trs homomorphisms.}
Suppose\sss that\sss $\mathcal{U}$\sss is\dss a star\sss finite proper covering of\dss $X$\nnsp.\oss
Let\sss $N_{\trf \mathcal{U}}$\sss be\sss the nerve of\dss $\mathcal{U}$\dnsp.\oss
Then\sss there\dss is\dss a canonical\qss 
\emph{Leray\trs homomorphism}\qss\vspace{1.5pt}
\[
\quad
l_{\dff \mathcal{U}}\dff \colon\dff
H_{\dff *}^\cf\trf(\trf X\trf)
\qff \ttoo\qff
H_{\dff *}^\iinf\dff(\trf N_{\trf \mathcal{U}} \trf)
\pff.
\]

\vspace{-10.5pt}
See\dss Section\qss \ref{cf-homology}.\oss
The\sss first\trs Leray\trs theorem\sss for $H_{\dff *}^\cf\trf(\trf X\trf)$\sss
is\dss the following\sss theorem.\oss

\myuppar{Theorem\qss A.}
\emph{Suppose\sss that\qss $\mathcal{U}$ is\dss a star\sss finite proper covering,\oss
and\dss that\qss $\mathcal{U}$ is\dss countable and\dss weakly $l_{\dff 1}$\dnsp-acyclic.\oss
If\qss a\sss compactly\dss finite\sss homology\sss class\sss
$h\qff \in\pff H_{\dff *}^\cf\trf(\trf X\trf)$\sss
belongs\sss to\sss the kernel\dss of\qss the\dss
Leray\trs homomorphism\dss $l_{\dff \mathcal{U}}$\nsp,\oss
then\dss $\norm{h}\off =\off 0$\nnsp.\oss}\vspace{6pt}

The second\trs Leray\trs theorem\dss is\dss concerned\sss with
almost\dss weakly $l_{\dff 1}$\dnsp-acyclic coverings.

\myuppar{Theorem\qss B.}
\emph{Suppose\sss that\qss $\mathcal{U}$ is\dss a star\sss finite proper covering,\oss
and\dss that\qss $\mathcal{U}$ is\dss countable and\dss
almost\dss weakly $l_{\dff 1}$\dnsp-acyclic.\oss
If\qss a\sss compactly\dss finite\sss homology\sss class\sss
$h\qff \in\pff H_{\dff *}^\cf\trf(\trf X\trf)$\sss
belongs\sss to\sss the kernel\dss of\qss the\dss
Leray\trs homomorphism\dss $l_{\dff \mathcal{U}}$\nsp,\oss
then\dss $\norm{h}\off <\off \infty$\nnsp.\oss}\vspace{6pt}

See\dss Theorems\qss \ref{norm-zero-weakly-acyclic}\qss
and\qss \ref{norm-finite-almost-weakly-acyclic}\qss
respectively.\oss
These\sss results are motivated\dss by\trs Gromov's\qss
\emph{Van\-ish\-ing-Finiteness\dss theorem}.\oss 
See\qss \cite{gro},\oss Section\qss 4.2.\oss 
Like\trs Leray\trs theorems of\qss \cite{i3},\pss
Theorems\qss A\qss and\qss B\qss are deduced\sss from an abstract\trs Leray\trs theorem,\oss
Theorem\qss \ref{e-acyclic-open}.\oss
The proofs are elementary\sss and are based on an adaptation of\trs the
methods of\pss \cite{i1},\qss \cite{i2},\qss and\qss \cite{i3}\qss to compactly\sss finite chains.\oss
The same methods\sss lead\sss to a proof\dss of\trs the\dss
Van\-ish\-ing-Finiteness\dss theorem.\oss
See\dss Section\qss \ref{extensions}.\oss

The assumptions of\qss Theorems\qss \ref{norm-zero-weakly-acyclic}\qss
and\qss \ref{norm-finite-almost-weakly-acyclic}\qss 
are much\sss weaker\dss than\dss Gromov's.\oss
In\dss the\dss Van\-ish\-ing-Finiteness\dss theorem\sss the space\sss $X$\sss
is\dss assumed\sss to be a manifold,\oss
instead of\dss $l_{\dff 1}$\dnsp-acyclicity\sss 
the stronger\sss amenability\dss property\dss is\dss used,\oss
and\dss it\dss is\dss assumed\dss that\sss
$h\qff \in\pff H_{\dff n}^\cf\trf(\trf X\trf)$\sss
for some $n$\sss strictly\dss larger\sss than\sss the dimension of\dss
$N_{\trf \mathcal{U}}$\nsp.\oss

Gromov's\trs proof\dss of\trs the\dss
Van\-ish\-ing-Finiteness\dss theorem\sss
was recently\sss reconstructed\dss by\dss
R.\dss Frigerio\dss and\trs M.\dss Moraschini\qss \cite{fm}.\oss
Their\sss proof\dss is\dss based on\trs Gromov's\trs theories of\dss
multicomplexes and of\dss diffusion of\dss chains,\oss 
and\dss is\dss far from\sss being elementary.\oss
Technical\sss difficulties forced\trs Frigerio\sss and\dss Moraschini\dss
to\sss consider only\dss triangulable spaces,\oss
although\dss they\sss conjectured\dss that\dss this assumption\dss is\dss
superfluous.\oss
Theorems\qss \ref{norm-zero-amenable-gromov}\qss and\qss
\ref{norm-finite-almost-compactly-amenable-gromov}\qss
together\sss imply\sss this conjecture.\oss

\myuppar{Removing\sss subspaces.}
Let\sss $X$\sss be a\sss topological\sss space,\oss
and\dss let\sss $Y\qff \subset\qff X$\sss be a closed subset.\oss
There exists natural\qss (in an\sss informal\sss sense)\qss chain\dss maps\vspace{3pt}
\[
\quad
r_{\dff \smallsetminus\dff Y}
\qff \colon\dff
C_{\dff \bullet}^\cf\dff(\trf X\trf)
\qff \ttoo\qff
C_{\dff \bullet}^\cf\dff(\trf X\qff \smallsetminus\qff Y\trf)
\pff.
\]

\vspace{-9pt}
See\dss Section\qss \ref{removing}.\oss
The maps\sss $r_{\dff \smallsetminus\dff Y}$\sss depend on many choices,\oss
but\dss the maps\vspace{3pt}
\[
\quad
\hry\dff \colon\dff
H_{\dff *}^\cf\dff(\trf X\trf)
\qff \ttoo\qff
H_{\dff *}^\cf\dff(\trf X\qff \smallsetminus\qff Y\trf)
\pff.
\]

\vspace{-8pt}
induced\sss by\sss $\ry$\sss does not\sss depend on\dss these choices.\oss
Suppose now\sss that\sss $Y$\sss is\dss presented as\sss the union of\dss a\sss family\sss
$\mathcal{Z}$\sss of\dss pair-wise disjoint\sss compact\sss subspaces of\dss $X$\nnsp.\oss
Suppose\sss further\dss that\dss for\sss every\sss
$Z\qff \in\qff \mathcal{Z}$\sss a compact\sss
neighborhood\dss $C_{\trf Z}$\dss of\dss $Z$\sss is\dss given,\oss
and\dss that\dss the neighborhoods\dss $C_{\trf Z}$\sss are pair-wise disjoint.\oss
Suppose\sss that\sss every\sss $C_{\trf Z}$\sss 
is\trs Hausdorff\trs and\dss path connected.\oss

\myuppar{Theorem\qss C.}
\emph{Suppose\sss that\sss $\mathcal{Z}$ is\dss countable
and\dss the\sss family\sss of\dss sets\sss $C_{\trf Z}$\sss
is\dss weakly\sss $l_{\dff 1}$\dnsp-acyclic.\oss 
Then\qss
$\norm{\ry\trf(\dff h\trf)}
\off \geq\off
\norm{h}$\pss
for every\dss homology\sss class\qss 
$h\qff \in\qff H_{\dff n}^\cf\dff(\trf X\trf)$\nnsp.\oss}\vspace{6pt}

See\dss Theorem\qss \ref{removing-lone}.\oss
Implicitly\dss this\sss theorem\dss is\dss concerned\sss with\dss the covering of\dss $X$\sss
by\dss $X\qff \smallsetminus\qff Y$\sss and\dss the sets\sss $C_{\trf Z}$\nsp.\oss
Since\sss this covering\dss is\dss very simple,\oss there\dss is\dss no need\dss to
involve\sss it\dss or\dss related\sss double complexes explicitly.\oss
Theorem\qss C\qss was motivated\dss by\trs Gromov's\qss
\emph{Cutting-of\qss theorem}\pss from\qss \cite{gro},\oss
Section\qss 4.2,\oss and easily\dss implies\sss its\sss first\sss claim.\oss
See\dss Section\qss \ref{removing}.\oss

\newpage
\mysection{A\pss Leray\qss theorem\qss for\pss infinite\qss chains}{homological-leray-infinite}

\myuppar{Generalized\sss chains.}
Let $X$ be a\sss topological\sss space.\oss
Let $\sub X$ be\sss the category\sss having\dss subspaces of\trs $X$\dss
as objects and\dss inclusions\dss
$Y\qff \subset\pff Z$\dss as morphisms.\oss
Let\dss $e_{\dff \bullet}$\dss be\dss a covariant\dss functor\dss from\dss
$\sub X$\dss to
augmented chain complexes of\dss modules over a ring\dss $R$\nnsp.\oss
The\sss functor\dss $e_{\dff \bullet}$\dss assigns\sss to a subspace $Z\qff \subset\pff X$ 
a\sss complex\vspace{-2pt}
\begin{equation}
\label{b-sigma}
\quad
\begin{tikzcd}[column sep=large, row sep=normal]\dis
0 \arrow[r]
& 
R  \arrow[l]
& 
e_{\trf 0}\dff(\trf Z \trf) 
\arrow[l, "\dis \off d_{\trf 0}\qff"']
& 
e_{\trf 1}\dff(\trf Z \trf)
\arrow[l, "\dis \off d_{\trf 1}\off"']
&   
e_{\trf 2}\dff(\trf Z \trf)
\arrow[l, "\dis \off d_{\trf 2}\off"']
&
\off \ldots \off\off,
\arrow[l, "\dis \off d_{\trf 3}\off"']
\end{tikzcd}
\end{equation}

\vspace{-9pt}
For every\dss $Y\qff \subset\qff Z$\dss there\dss is\dss a\qss
\emph{inclusion\dss morphism}\qss
$e_{\trf \bullet}\dff(\trf Z \trf)
\dff \ttoo\dff
e_{\trf \bullet}\dff(\trf Y \trf)$\dnsp.\oss
Elements\sss of\dss $e_{\fff q}\dff(\trf Z \trf)$
are\sss thought\sss as\qss \emph{generalized $q$\dnsp-chains}\qss of\dss $Z$\nnsp.\oss

\myuppar{The double complex of\dss a covering.}
Let $\mathcal{U}$ be a star\sss finite covering\sss of\dss $X$\sss and\dss let\sss
$N$\sss be\sss its\sss nerve.\oss
For\sss $p\qff \geq\qff 0$\sss let\dss $N_{\dff p}$\dss 
be\sss the set\sss of\dss $p$\dnsp-dimensional\sss simplices of\dss $N$\nnsp.\oss
For\sss $p\fff,\pff q\qff \geq\qff 0$\sss let\vspace{3pt}\vspace{-1pt}
\[
\quad
c_{\dff p}\dff(\trf N\fff,\pff e_{\dff q} \trf)
\off\off =\off\off
\prod\nolimits_{\qff \sigma\qff \in\qff N_{\dff p}}\qff 
e_{\dff q}\fff\left(\trf \num{\sigma} \trf\right)
\pff.
\]

\vspace{-9pt}\vspace{-1pt}
So,\oss an element 
$c\qff \in\qff
c_{\dff p}\dff(\trf N\fff,\pff e_{\dff q} \trf)$
is\dss a\sss family of\qss \emph{generalized
$q$\dnsp-chains}\vspace{3pt}\vspace{-1pt}
\begin{equation*}
\quad
c\dff \colon\dff
\sigma
\off \longmapsto\off
c_{\dff \sigma}
\off \in\off
e_{\dff q}\fff\left(\trf \num{\sigma} \trf\right)
\pff,
\end{equation*}

\vspace{-9pt}\vspace{-1pt}
where\dss $\sigma\qff \in\qff N_{\dff p}$\nsp,\oss
thought\sss as an\dss infinite formal\sss sum\vspace{3pt}\vspace{-1pt}
\[
\quad 
c
\off =\off
\sum\nolimits_{\qff \sigma\qff \in\qff N_{\dff p}}\qff
c_{\dff \sigma}
\pff.
\]

\vspace{-9pt}\vspace{-1pt}
For every\dss $p\qff >\qff 0$\dss 
(and\sss sometimes for $p\off =\off 0$ also)\qss
there\dss is\dss a canonical\dss morphism\qss\vspace{3pt}\vspace{-1pt}
\[
\quad 
\delta_{\fff p}\dff \colon\dff
c_{\dff p}\dff(\dff N\fff,\pff e_{\dff \bullet} \dff)
\qff \ttoo\qff
c_{\dff p\dff -\dff 1}\dff(\dff N\fff,\pff e_{\dff \bullet} \dff)
\pff,
\]

\vspace{-9pt}\vspace{-1pt}
defined as follows.\oss
Let\dss $\sigma\qff \in\qff N_{\dff p}$\nsp.\oss
For each\dss face\dss $\partial_{\dff i}\dff\sigma$\dss
there\dss is\dss an\dss inclusion\dss morphism\dss\vspace{3pt}\vspace{-1pt}
\[
\quad
\Delta_{\qff \sigma\fff,\dff i}
\pff \colon\pff
e_{\dff \bullet}\dff(\trf \num{\sigma} \trf)
\qff \ttoo\qff
e_{\dff \bullet}\dff(\trf \num{\partial_{\dff i}\dff\sigma} \trf)
\pff. 
\]

\vspace{-9pt}\vspace{-1pt}
For\dss 
$c_{\dff \sigma}
\qff \in\qff 
e_{\dff q}\fff\left(\trf \num{\sigma} \trf\right)$\dss
we set\vspace{3pt}\vspace{0.15pt}
\[
\quad
\delta_{\fff p}\dff \left(\trf c_{\dff \sigma}\trf\right)
\qff \colon\qff
\partial_{\dff i}\dff\sigma
\off \longmapsto\off
(\dff -\qff 1 \dff)^{\dff i}\qff
\Delta_{\qff \sigma\fff,\dff i}\dff\left(\trf
c_{\dff \sigma}
\trf\right)
\off \in\off
e_{\dff q}\fff\left(\trf \num{\partial_{\dff i}\dff\sigma} \trf\right)
\quad\
\mbox{and}
\]

\vspace{-36pt}\vspace{0.15pt}
\[
\quad
\delta_{\fff p}\dff \left(\trf c_{\dff \sigma}\trf\right)
\qff \colon\qff
\tau
\off \longmapsto\off
0
\]

\vspace{-9pt}\vspace{0.15pt}
if\qss 
$\tau\off \neq\off \partial_{\dff i}\dff\sigma$\qss
for every\sss $i$\nnsp.\oss
The map\dss $\delta_{\fff p}$\dss extends\sss  
to\sss the direct\sss product\dss
$c_{\dff p}\dff(\trf N\fff,\pff e_{\dff q} \trf)$\dss
by\dss linearity\halfff.\oss
In order\dss to see\sss that\sss such an extension\dss to be well\sss defined\sss
we need\dss to know\dss that\dss for every\sss
$\rho\qff \in\pff N_{\dff p\dff -\dff 1}$\sss
only\sss a\sss finite number of\dss expressions\sss
$(\dff -\qff 1 \dff)^{\dff i}\qff
\Delta_{\qff \sigma\fff,\dff i}\dff\left(\trf
c_{\dff \sigma}
\trf\right)$\sss
need\dss to be summed\dss
to get\dss the value of\dss
$\delta_{\fff p}\dff \left(\trf c_{\dff \sigma}\trf\right)$\sss
on\sss $\rho$\nnsp.\oss
But 
$(\dff -\qff 1 \dff)^{\dff i}\qff
\Delta_{\qff \sigma\fff,\dff i}\dff\left(\trf
c_{\dff \sigma}
\trf\right)$
enters\sss this sum only\trs if\qss
$\rho\off =\off \partial_{\dff i}\dff\sigma$\sss
for\sss some\dss
$\sigma\qff \in\pff N_{\dff p}$\sss and\sss some $i$\nnsp.\oss
Since\sss the covering\sss $\mathcal{U}$\sss is\dss star\sss finite,\oss
its\sss nerve\sss $N$\sss is\dss also star\sss finite,\oss
and\dss hence\sss there\dss is\dss indeed only\sss a finite number of\dss
such\sss $\rho\fff,\qff i$\nnsp.\oss
Therefore\sss $\delta_{\fff p}$\sss
is\dss indeed\sss well\sss defined.\oss
As usual,\oss we agree\sss that\sss $\num{\varnothing}\off =\off X$\nnsp,\oss
but\dss this argument\sss does not\sss work for\sss
$p\off =\off 0$\sss and\sss $\rho\off =\off \varnothing$\nnsp.\oss
In\sss order\sss to define\sss $\delta_{\dff 0}$\sss
one needs\sss to be able\sss to speak about\dss infinite generalized chains.\oss

The fact\dss that\sss 
each\dss 
$\Delta_{\qff \sigma\fff,\dff i}$\dss 
are\dss
morphisms of\dss complexes implies\sss that\dss
$\delta_{\fff p}$\dss is\dss a\sss morphism also of\dss complexes.\oss
The\qss \emph{double complex}\pss
$c_{\trf \bullet}\fff(\trf N\fff,\pff e_{\trf \bullet} \trf)$\qss 
\emph{of\qss the covering}\qss $\mathcal{U}$\qss
is\dss the double complex\vspace{4.25pt}
\begin{equation}
\label{double-complex-covering-homology}
\qquad
\begin{tikzcd}[column sep=large, row sep=huge]\dis
c_{\trf 0}\dff(\trf N\fff,\pff e_{\trf 0} \trf) 
&   
c_{\trf 0}\dff(\trf N\fff,\pff e_{\trf 1} \trf) \arrow[l]
&
\off \ldots \off \arrow[l]
\\
c_{\dff 1}\dff(\trf N\fff,\pff e_{\trf 0} \trf) 
\arrow[u]
&   
c_{\dff 1}\dff(\trf N\fff,\pff e_{\trf 1} \trf) \arrow[l]
\arrow[u]
&
\off \ldots \off \arrow[l]
\\
c_{\trf 2}\dff(\trf N\fff,\pff e_{\trf 0} \trf) 
\arrow[u]
&   
c_{\trf 2}\fff(\trf N\fff,\pff e_{\trf 1} \trf) \arrow[l]
\arrow[u]
&
\off \ldots \off \arrow[l]
\\
\ldots \vphantom{C_{\trf 2}\fff(\trf N\fff,\pff e_{\trf 1} \trf)} 
\arrow[u] 
&
\ldots \vphantom{C_{\trf 2}\fff(\trf N\fff,\pff e_{\trf 1} \trf)} 
\arrow[u]
&
\quad \vphantom{C_{\trf 2}\fff(\trf N\fff,\pff e_{\trf 1} \trf)} 
,
\end{tikzcd}
\end{equation}

\vspace{-10.75pt}
where\sss the horizontal\sss arrows are\sss the products of\trs the maps\dss $d_{\dff i}$\dss
and\dss the vertical\sss arrows are\sss the maps\dss $\delta_{\dff i}$\nsp.\oss
Let\trs 
$t_{\trf \bullet}\dff(\trf N\fff,\pff e \trf)$\dss
be\sss the\sss total\sss complex of\dss
$c_{\trf \bullet}\fff(\trf N\fff,\pff e_{\trf \bullet} \trf)$\nnsp.\oss
Let\dss 
$C_{\dff \bullet}^\iinf\fff(\trf N\trf)
\off =\off
C_{\dff \bullet}^\iinf\dff(\trf S\fff,\qff R\trf)$\dss
be\sss the complex of\qss \emph{infinite}\qss 
simplicial\sss chains of\trs $N$\dss with coefficients\sss
in\dss $R$\nnsp.\oss
It\dss is\dss well\sss defined\dss 
because\sss $N$\sss is\dss star\sss finite.\oss
Let\dss 
$H_{\dff \bullet}^\iinf\dff(\trf N\trf)
\off =\off
H_{\dff \bullet}^\iinf\dff(\trf N\fff,\qff R\trf)$\sss
be\sss the homology\sss of\trs this complex.\oss
Since\sss $\delta_{\dff 0}$\sss is\dss not\sss defined,\oss
we will\sss replace\sss $e_{\dff \bullet}\dff(\trf X\trf)$\sss
by\dss the cokernel\sss
$e_{\dff \bullet}^\lf\dff(\qff X\fff,\pff \mathcal{U} \trf)$\sss
of\trs the homomorphism\vspace{3pt}
\[
\quad
\delta_{\dff 1}\dff \colon\dff
c_{\dff 1}\dff(\dff N\fff,\pff e_{\dff \bullet} \dff)
\qff \ttoo\qff
c_{\trf 0}\dff(\dff N\fff,\pff e_{\dff \bullet} \dff)
\pff.
\]

\vspace{-9pt}
Then\sss we can define\sss $\delta_{\dff 0}$\sss as\sss the canonical\dss map\sss
$c_{\trf 0}\dff(\dff N\fff,\pff e_{\dff \bullet} \dff)
\qff \ttoo\qff
e_{\dff \bullet}^\lf\dff(\qff X\fff,\pff \mathcal{U} \trf)$\nnsp.\oss
Since\qss (\ref{double-complex-covering-homology})\qss is\dss commutative,\oss
the maps\sss $d_{\dff i}$\sss induce canonical\dss maps\vspace{3pt}
\[
\quad
d_{\dff i}\dff \colon\dff
e_{\dff i}^\lf\dff(\qff X\fff,\pff \mathcal{U} \trf)
\qff \ttoo\qff
e_{\dff i\dff -\dff 1}^\lf\dff(\qff X\fff,\pff \mathcal{U} \trf)
\]

\vspace{-9pt}
turning\sss
$e_{\dff \bullet}^\lf\dff(\qff X\fff,\pff \mathcal{U} \trf)$\sss
into a complex.\oss
Let\sss
$\widetilde{H}_{\dff \bullet}^\lf\dff (\qff X\fff,\pff \mathcal{U} \trf)$\sss
be\sss the homology\sss of\trs this complex.\oss
The boundary\dss maps\dss $d_{\trf 0}$\dss and\dss $\delta_{\dff 1}$\dss 
lead\dss to morphisms\dss\vspace{3pt}
\[
\quad
\lambda_{\dff e}\qff \colon\qff
t_{\dff \bullet}\dff(\trf N\fff,\pff e \trf) 
\qff \ttoo\qff
C_{\dff \bullet}^\iinf\dff(\trf N \trf)
\quad\
\mbox{and}\quad\
\tau_{\dff e}\qff \colon\qff
t_{\dff \bullet}\dff(\trf N\fff,\pff e \trf) 
\qff \ttoo\qff
e_{\dff \bullet}^\lf\dff(\qff X\fff,\pff \mathcal{U} \trf)
\pff,
\]

\vspace{-9pt}
where\dss it\dss is\dss understood\dss that\dss the augmentation\dss term\dss 
is\dss  removed\dss from\dss 
$e_{\dff \bullet}^\lf\dff(\qff X\fff,\pff \mathcal{U} \trf)$\nnsp.\oss

\myuppar{Acyclic coverings.}
Clearly,\pss $\mathcal{U}^{\dff \cap}$\sss is\dss the collection of\dss
all\sss sets of\trs the form\dss $\num{\sigma}$\sss
with\sss $\sigma\off \neq\off \varnothing$\nnsp.\oss
The covering\dss $\mathcal{U}$\dss is\dss said\dss to be\dss 
\emph{$e_{\fff \bullet}$\nsp\dnsp-acyclic}\oss if\pss
$e_{\dff \bullet}\dff(\trf Z\trf)$\dss is\dss exact\dss for every\dss 
$Z\qff \in\qff \mathcal{U}^{\dff \cap}$\dnsp.\oss

\mypar{Lemma.}{acyclic-coverings}
\emph{If\pss $\mathcal{U}$\dss is\dss star\dss finite and\qss 
$e_{\fff \bullet}$\nsp\dnsp-acyclic,\oss 
then\qss
$\lambda_{\dff e}\dff \colon\dff
t_{\dff \bullet}\dff(\trf N\fff,\pff e \trf) 
\qff \ttoo\qff
C_{\dff \bullet}^\iinf\dff(\trf N \trf)$\qss
induces an\dss isomorphism of\qss homology\dss groups.\oss}

\proof
If\qss $\mathcal{U}$\dss is\pss 
$e_{\fff \bullet}$\nsp\dnsp-acyclic,\oss
then\dss for every\dss  
simplex\dss $\sigma\off \neq\off \varnothing$\dss the complex\qss
(\ref{b-sigma})\qss is\dss exact\halfff.\oss
Since\sss the\sss term-wise products of\dss exact\sss sequences are exact,\oss
this implies\sss that\dss every\dss row\sss of\trs the double complex\qss
(\ref{double-complex-covering-homology})\qss is\dss exact\sss
and\dss $d_{\trf 0}$\dss 
induces an\dss isomorphism of\trs the complex\dss
$C_{\dff \bullet}^\iinf\dff(\trf N \trf)$\dss 
with\dss the kernel\dss of\trs the morphism of\dss complexes\dss 
$d_{\dff 1}\dff \colon\dff
c_{\dff \bullet}\dff(\trf N\fff,\pff e_{\dff 1} \trf)
\qff \ttoo\qff
c_{\dff \bullet}\dff(\trf N\fff,\pff e_{\trf 0} \trf)$\nnsp.\oss
It\dss remains\sss to apply\sss a well\dss known\dss 
theorem about\sss double complexes.\oss
See\trs Theorem\qss A.2\qss in\qss \cite{i3}.\oss  \eproof

\myuppar{Infinite singular\sss chains.}
Suppose\sss that\dss a\sss space $\Delta$ is\dss fixed and\dss
maps\dss $s\dff \colon\dff \Delta\qff \ttoo\qff Y$\dss
are\sss treated as singular simplices.\oss
A\qss \emph{finite singular\dss chain}\pss is\dss a finite formal\sss sum
of\dss singular\sss simplices with coefficients in\dss $R$\nnsp.\oss
The $R$\dnsp-module of\dss
finite singular chains\dss
in\sss $Y$\sss is\dss denoted\dss by\sss $c\trf(\trf Y\trf)$\nnsp.\oss
An\qss \emph{infinite singular\dss chain}\pss is\dss a finite or\sss infinite formal\sss sum
of\dss singular\sss simplices with coefficients in\dss $R$\nnsp,\oss
and\dss the $R$\dnsp-module of\dss
infinite singular chains\dss
in\sss $Y$\sss is\dss denoted\dss by\sss $c^\iinf\trf(\trf Y\trf)$\nnsp.\oss

Let\dss us\sss turn\sss to singular chains\sss in\sss $X$\sss and subsets of\dss $X$\nnsp.\oss
A singular\sss simplex $s$\sss in\sss $X$\sss is\dss called\qss \emph{small}\oss if\dss 
$s\dff(\trf \Delta\trf)\qff \subset\qff U$\sss
for\sss some\sss $U\qff \in\qff \mathcal{U}$\nnsp,\oss
and an\dss infinite singular chain\sss in\sss $X$\sss
is\dss called\qss \emph{small}\oss if\dss all\sss its\sss singular\sss simplices
with non-zero coefficients are small.\oss
Suppose\sss that\dss for every\sss $U\qff \in\qff \mathcal{U}$\sss
a\sss finite singular chain\sss
$\gamma_{\dff U}\qff \in\qff c\trf(\trf U\trf)$\sss
is\dss given.\oss
Then,\oss since\sss $\mathcal{U}$\dss is\dss
a star\sss finite,\oss the sum\vspace{3pt}\vspace{0.5pt}
\[
\quad
\gamma
\off =\off
\sum\nolimits_{\qff U\qff \in\qff \mathcal{U}}\qff
\gamma_{\dff U}
\]

\vspace{-8.25pt}
is\dss a\sss well\sss defined\dss infinite chain.\oss
Clearly,\pss $\gamma$\sss is\dss small.\oss 
If\dss $\mathcal{U}$\sss is\dss an\sss open covering\halfff,\oss
then\sss $\gamma$\sss is\dss locally\dss finite.\oss
When\dss a chain\sss $\gamma$\sss can\sss be represented\dss by\sss such sum,\oss
we\sss say\dss that\sss $\gamma$\sss is\dss \emph{$\mathcal{U}$\nsp\dnsp-finite}.\oss
Let\dss $c^\lf\trf(\qff X\fff,\qff \mathcal{U} \trf)$\sss 
be\sss the $R$\dnsp-module of\dss $\mathcal{U}$\nsp\dnsp-finite chains.\oss
Let\dss us\sss consider\sss now\sss the modules\vspace{1.5pt}\vspace{1.25pt}
\[
\quad
c_{\dff p}\dff(\trf N\fff,\pff c \trf)
\off\off =\off\off
\prod\nolimits_{\qff \sigma\qff \in\qff N_{\dff p}}\qff 
c\dff\left(\trf \num{\sigma} \trf\right)
\pff,
\]

\vspace{-10.5pt}\vspace{1.25pt}
where\dss $p\qff \geq\qff -\qff 1$\nnsp.\oss
The maps\qss 
$\delta_{\fff p}\dff \colon\dff
c_{\dff p}\dff(\trf N\fff,\pff c \trf)
\qff \ttoo\qff
c_{\dff p\dff -\dff 1}\dff(\trf N\fff,\pff c \trf)$\nnsp,\pss
$p\qff >\qff 0$\nnsp,\oss
are defined as before.\oss
Moreover\halfff,\oss
now\sss we can define\sss $\delta_{\trf 0}$\sss in\dss the same way,\oss
except\dss that\dss now\dss the\sss target\sss of\dss $\delta_{\trf 0}$\sss
is\sss $c^\iinf\dff(\trf X \trf)$\nnsp,\oss
not\sss $c\trf(\trf X\trf)$\nnsp.\oss
Clearly,\pss $c^\lf\trf(\qff X\fff,\pff \mathcal{U} \trf)$\sss
is\dss equal\dss to\sss the image of\vspace{1.5pt}\vspace{1.25pt}
\[
\quad
\delta_{\trf 0}\dff \colon\dff
c_{\trf 0}\dff(\dff N\fff,\pff c \dff)
\qff \ttoo\qff
c^\iinf\dff(\trf X \trf)
\pff,
\]

\vspace{-10.5pt}\vspace{1.25pt}
and\sss $\delta_{\trf 0}$\sss is\dss the composition of\trs the inclusion\sss
$c^\lf\trf(\qff X\fff,\pff \mathcal{U} \trf)
\qff \ttoo\qff
c^\iinf\dff(\trf X \trf)$\sss
with a\sss canonical\dss map\vspace{1.5pt}\vspace{1.25pt}
\[
\quad
\overline{\delta}_{\trf 0}\dff \colon\dff
c_{\trf 0}\dff(\trf N\fff,\pff c \trf)
\qff \ttoo\qff
c^\lf\dff(\qff X\fff,\pff \mathcal{U} \trf)
\pff.
\]

\vspace{-10.5pt}\vspace{1.25pt}
Clearly,\oss $\overline{\delta}_{\trf 0}$\sss is\dss surjective.\oss
As usual,\pss
$\delta_{\fff p\dff -\dff 1}\dff \circ\qff \delta_{\fff p}
\off =\off
0$\dss
for\sss every\dss $p\qff \geq\qff 1$\nnsp.\oss

\mypar{Lemma.}{homology-columns-are-exact}
\emph{If\pss $\mathcal{U}$\dss is\dss star\dss finite,\oss
then\dss the\qss following\dss sequence\dss is\dss exact\dff:}\vspace{-1.5pt}
\[
\quad
\begin{tikzcd}[column sep=large, row sep=normal]\dis
0
&
c^\lf\dff(\trf X\fff,\pff \mathcal{U}\trf) \arrow[l]
&
c_{\trf 0}\dff(\trf N\fff,\pff c \trf) \arrow[l, "\dis \off \overline{\delta}_{\trf 0}\qff"']
& 
c_{\trf 1}\dff(\trf N\fff,\pff c \trf) \arrow[l, "\dis \off \delta_{\trf 1}\qff"'] 
&   
\off \ldots\off . \arrow[l, "\dis \off \delta_{\trf 2}"'] 
\end{tikzcd}
\]

\vspace{-9pt}
\proof
It\dss is\dss sufficient\dss to\sss construct\sss 
a contracting\sss chain\dss homotopy\vspace{4.5pt}
\[
\quad
k_{\qff 0}\dff \colon\dff
c^\lf\dff(\trf X\fff,\pff \mathcal{U}\trf) 
\qff \ttoo\qff
c_{\trf 0}\dff(\trf N\fff,\pff c \trf)\dff,\quad\ 
\]

\vspace{-34.5pt}
\[
\quad
k_{\trf p}\dff \colon\dff
c_{\trf p}\dff(\trf N\fff,\pff c \trf) 
\qff \ttoo\qff
c_{\trf p\dff +\dff 1}\dff(\trf N\fff,\pff c \trf) 
\dff, 
\]

\vspace{-7.5pt}
where\dss $p\qff \geq\qff 0$\nnsp.\oss
The construction\dss is\dss almost\dss the same as in\sss the case of\dss
direct\sss sums\qss (instead of\dss products).\oss
Cf.\pss \cite{i3},\oss Lemma\qss 3.1.\oss
For every\sss small\sss singular simplex 
$s$
let\dss us\sss choose a subset\dss $U_{\dff s}\qff \in\qff \mathcal{U}$\dss
such\dss that\dss 
$s\trf(\dff \Delta\dff) 
\qff \subset\pff 
U_{\dff s}$\dss
and\dss let\sss $u_{\fff s}$\sss be\sss the corresponding\dss vertex of\sss $N$\nnsp.\oss
If\trs 
$s\dff(\dff \Delta\trf)
\qff \subset\qff
\num{\sigma}$\dss
for some\dss $\sigma\qff \in\qff N_{\dff p}$\nsp,\oss
then\sss $s\dff *\dff \sigma$\sss 
denotes $s$
considered as an element\dss of\trs $c\trf(\trf \num{\sigma}\trf)$\nnsp.\oss

Let\dss us\sss define\sss 
$k_{\trf p}$
on\dss the chains of\trs the form\sss
$s\dff *\dff \sigma$\sss first.\oss 
Suppose\sss that\dss
$\sigma\qff \in\qff N_{\dff p}$
and\sss $s$\dss be a singular $q$\dnsp-simplex
such\dss that\dss
$s\dff(\dff \Delta^{\fff q}\trf)
\qff \subset\qff
\num{\sigma}$\nnsp.\oss
Let\dss
$\rho
\off =\off
\sigma\dff \cup\dff \{\trf u_{\dff s}\trf\}$\nnsp.\oss
Then\dss
$s\dff(\dff \Delta\trf)
\qff \subset\qff
\num{\sigma}\qff \cap\qff U_{\dff s}
\off =\off
\num{\rho}$\dss
and,\oss in\dss \particular\halfff,\pss $\rho$\dss is\dss a simplex.\oss
If\trs $u_{\dff s}\qff \in\pff \sigma$\nnsp,\oss
then\dss $\rho$\dss is\dss a $p$\dnsp-simplex.\oss
Otherwise,\pss $\rho$\dss is\dss a $(\dff p\qff +\qff 1\dff)$\dnsp-simplex and\dss
$\sigma\off =\off \partial_{\fff a}\trf \rho$\dss for some\dss $a$\nnsp.\oss 
Let\vspace{4.5pt}
\[
\quad
k_{\trf p}\dff(\dff s\dff *\dff \sigma\trf)
\off =\off
0
\quad\
\mbox{if}\quad\
u_{\dff s}\qff \in\pff \sigma
\pff,
\]

\vspace{-34.5pt}
\[
\quad
k_{\trf p}\dff(\dff s\dff *\dff \sigma\trf)
\off =\off
(\qff -\qff 1\trf)^{\dff a}\qff 
s\dff *\dff \rho\off \in\off c\trf(\trf \num{\rho}\trf)
\quad\
\mbox{if}\quad\
u_{\dff s}\qff \not\in\pff \sigma
\pff.
\]

\vspace{-7.5pt}
As\sss in\sss the case of\dss $\delta_{\dff p}$\nsp,\oss
the star\sss finitness of\dss $\mathcal{U}$ and $N$ 
allows\sss to extend\sss $k_{\trf p}$\sss 
to $c_{\trf p}\dff(\trf N\fff,\pff c \trf)$\sss
by\dss linearity.\oss
In\sss order\sss to verify\dss that $k_{\dff \bullet}$\sss is\dss
a contracting\sss homotopy\dss it\dss is\dss 
sufficient\dss to check\dss that\vspace{4pt}
\[
\quad
\delta_{\dff p\dff +\dff 1}\trf \bigl(\trf k_{\dff p}\dff(\trf \gamma\trf)
\qff\bigr)
\off +\off
k_{\dff p\dff -\dff 1}\trf \bigl(\qff \delta_{\dff p}\dff(\trf \gamma\trf)
\qff\bigr)
\off =\off
\gamma
\]

\vspace{-8pt}
when\sss $\gamma$\sss has\sss the\sss form\sss
$\gamma\off =\off s\dff *\dff \sigma$\nnsp.\oss
But\dss this case\dss is\dss exactly\dss the same
as for direct\sss sums.\oss  \eproof

\myuppar{Classical\sss singular\sss chains.}
The above discussion applies,\oss in\dss particular\halfff,\oss
to\sss the case\sss $\Delta\off =\off \Delta^q$\dnsp,\oss
the standard\sss geometric $q$\dnsp-simplex.\oss
In\dss this case we will\sss denote\sss
$c_{\dff p}\dff(\trf N\fff,\pff c \trf)$\sss
and\sss
$c^\lf\trf(\qff X\fff,\pff \mathcal{U} \trf)$\sss by\vspace{3pt}
\[
\quad
C_{\dff p}^\iinf\dff(\trf N\fff,\pff C_{\dff q} \trf)
\quad\
\mbox{and}\quad\
C_{\dff q}^\lf\trf(\qff X\fff,\pff \mathcal{U} \trf)
\]

\vspace{-9pt}
respectively.\oss
The boundary\sss maps\sss $d_{\dff i}$\sss turn\sss
$C_{\dff \bullet}^\lf\trf(\qff X\fff,\pff \mathcal{U} \trf)$\sss
into a complex.\oss
Let\dss $H_{\dff *}^\lf\dff(\qff X\fff,\pff \mathcal{U} \trf)$\sss
be\sss the homology\sss of\trs this complex.\oss
The morphisms\vspace{-1.5pt}
\[
\quad
\begin{tikzcd}[column sep=large, row sep=huge]\dis
C_{\dff \bullet}^\iinf\dff(\trf N \trf) 
& 
t_{\dff \bullet}\dff(\trf N\fff,\pff C \dff)
\arrow[l, "\dis \off\qff \lambda_{\qff C}"']
\arrow[r, "\dis \tau_{\trf C}\off"]
&
C_{\dff \bullet}^\lf\trf(\trf X\fff,\pff \mathcal{U}\trf)
\end{tikzcd}
\]

\vspace{-9pt}
lead\dss to homomorphisms of\dss cohomology\dss groups,\oss\vspace{-4pt}
\[
\quad
\begin{tikzcd}[column sep=large, row sep=huge]\dis
H_{\dff *}^\iinf\dff(\trf N \trf) 
& 
H_{\dff *}\dff(\trf N\fff,\pff C \dff)
\arrow[l, "\dis \off\qff \lambda_{\qff C\dff *}"']
\arrow[r, "\dis \tau_{\trf C\dff *}\off"]
&
H_{\dff *}^\lf\trf(\trf X\fff,\pff \mathcal{U}\trf)
\end{tikzcd}
\]

\vspace{-10pt}
where\dss 
$H_{\dff *}\dff(\trf N\fff,\pff C \dff)$\dss
is\sss the homology\sss of\trs
$t_{\dff \bullet}\dff(\trf N\fff,\pff C \dff)$\nnsp.\oss
If\dss $\mathcal{U}$\sss is\dss star\sss finite,\oss
Lemma\qss \ref{homology-columns-are-exact}\qss implies\sss that\dss
the columns of\qss (\ref{double-complex-covering-homology})\qss are exact.\oss
Together\sss with\sss the\sss already\sss used\dss theorem about\sss double complexes\sss
this implies\sss that\sss $\tau_{\trf C\dff *}$\sss is\dss an\sss isomorphism.\oss
This\sss leads\sss to\sss the canonical\dss homomorphism\vspace{3pt}
\begin{equation}
\label{leray-map}
\quad
\lambda_{\qff C\dff *}\dff \circ\dff \tau_{\trf C\dff *}^{\dff -\dff 1}
\qff \colon\qff
H_{\dff *}^\lf\trf(\trf X\fff,\pff \mathcal{U}\trf)
\off \ttoo\off
H_{\dff *}^\iinf\dff(\trf N \trf) 
\pff.
\end{equation}

\vspace{-9pt}
In\sss general,\oss one cannot\dss replace here\sss
$H_{\dff \bullet}^\lf\trf(\trf X\fff,\pff \mathcal{U}\trf)$\sss
by\sss some homology\sss independent\sss of\dss $\mathcal{U}$\dnsp.\oss

\myuppar{Comparing\sss classical\sss and\dss generalized chains.}
Suppose\sss that\dss the functor\sss $e_{\dff \bullet}$\sss is\dss
equipped\sss with a natural\dss transformation\sss
$\varphi_{\dff \bullet}\dff \colon\dff
C_{\dff \bullet}\qff \ttoo\qff e_{\dff \bullet}$\nnsp.\oss
Then\sss $\varphi$\sss induces a map\sss
$H_{\dff *}\dff(\trf Y\trf)
\qff \ttoo\qff
\widetilde{H}_{\dff *}\dff(\trf Y\trf)$\sss
for every\sss $Y\qff \in\qff \mathcal{U}$\nnsp,\oss
where\sss $\widetilde{H}_{\dff *}\dff(\trf Y\trf)$\sss
is\dss the\sss homology\sss of\trs the complex\sss
$e_{\dff \bullet}\dff(\trf Y\trf)$\nnsp.\oss
In\sss particular\halfff,\pss $\varphi$\sss
induces a map\sss
$H_{\dff *}\dff(\trf X\trf)
\qff \ttoo\qff
\widetilde{H}_{\dff *}\dff(\trf X\trf)$\nnsp,\oss
but\dss this\dss is\dss not\dss what\sss we are interested\sss in\dss now.\oss
Lemma\qss \ref{homology-columns-are-exact}\qss implies\sss that\dss
the complex\sss $C_{\dff \bullet}^\lf\trf(\trf X\fff,\pff \mathcal{U}\trf)$\sss
is\dss canonically\dss isomorphic\sss to\sss the cokernel\sss of\vspace{3pt}
\[
\quad
\delta_{\dff 1}\dff \colon\dff
C_{\trf 1}^\iinf\dff(\trf N\fff,\pff C_{\dff \bullet} \trf)
\qff \ttoo\qff
C_{\trf 0}^\iinf\dff(\trf N\fff,\pff C_{\dff \bullet} \trf)
\pff.
\]

\vspace{-9pt}
In\sss view of\trs the definition of\dss
$e_{\dff \bullet}^\lf\dff(\qff X\fff,\pff \mathcal{U} \trf)$\sss
this\sss leads\sss to a canonical\dss homomorphism\vspace{3pt}
\[
\quad
\varphi_{\dff \bullet}\dff \colon\dff
C_{\dff \bullet}^\lf\trf(\trf X\fff,\pff \mathcal{U}\trf)
\qff \ttoo\qff
e_{\dff \bullet}^\lf\dff(\qff X\fff,\pff \mathcal{U} \trf)
\]

\vspace{-9pt}
and\dss hence\sss to a\qss \emph{comparison\dss homomorphism}\dss\vspace{3pt}
\begin{equation}
\label{comparison-map}
\quad
\varphi_{\dff *}\dff \colon\dff
H_{\dff *}^\lf\trf(\trf X\fff,\pff \mathcal{U}\trf)
\qff \ttoo\qff
\widetilde{H}_{\dff *}^\lf\dff(\qff X\fff,\pff \mathcal{U} \trf)
\end{equation}

\vspace{-9pt}
in\sss homology\sss groups.\oss
This\dss is\dss the map\sss we are interested\dss in.\oss

\mypar{Theorem.}{e-acyclic-open}
\emph{If\pss $\mathcal{U}$\sss is\dss a\sss star\dss finite\dss 
$e_{\bullet}$\dnsp-acyclic\sss covering,\oss
then\dss the\sss comparison\dss homomorphism\qss
\textup{(\ref{comparison-map})}\qss
can\sss be\trs factored\dss through\dss 
the\sss canonical\trs 
homomorphism\qss
\textup{(\ref{leray-map})}.\oss}

\proof
The natural\dss transformation\sss $\varphi_{\dff \bullet}$\sss 
defines\sss homomorphisms\qss\vspace{1.5pt}
\[
\quad
\varphi_{\dff q}\dff\left(\trf \num{\sigma} \trf\right)\dff \colon\dff
C_{\dff q}\fff\left(\trf \num{\sigma} \trf\right)
\qff \ttoo\qff
e_{\dff q}\fff\left(\trf \num{\sigma} \trf\right)
\pff,
\]

\vspace{-10.5pt}
which,\oss in\sss turn,\oss lead\dss to a morphism\vspace{1.5pt}
\[
\quad
\varphi_{\dff \bullet\dff \bullet}\dff \colon\dff
C_{\dff \bullet}\dff(\trf N\fff,\pff C_{\dff \bullet} \trf)
\off \ttoo\off
C_{\dff \bullet}\dff(\trf N\fff,\pff e_{\dff \bullet} \trf)
\pff
\]

\vspace{-10.5pt}
of\dss double complexes.\oss 
In\dss turn,\pss 
$\varphi_{\dff \bullet\dff \bullet}$\dss leads\sss to a\sss morphism\dss
$\Phi_{\dff \bullet}\dff \colon\dff
C_{\dff \bullet}\dff(\trf N\fff,\pff C \trf)
\off \ttoo\off
t_{\dff \bullet}\dff(\trf N\fff,\pff e \trf)$\dss
of\trs total\sss complexes.\oss
Clearly\halfff,\oss the diagram\vspace{4.5pt}
\[
\quad
\begin{tikzcd}[column sep=boommm, row sep=boommm]\dis
C_{\dff \bullet}^\iinf\dff(\trf N \trf) 
\arrow[d, "\dis \qff ="]
& 
t_{\dff \bullet}\dff(\trf N\fff,\pff C \trf)
\arrow[d, "\dis \qff \Phi_{\dff \bullet}"]
\arrow[r, "\dis \off\qff \tau_{\trf C}\off"]
\arrow[l, "\dis \off\qff \lambda_{\trf C}"']
&
C_{\dff \bullet}^\lf\dff(\trf X\fff,\pff \mathcal{U}\trf)
\arrow[d, "\dis \qff \varphi_{\dff \bullet}"]
\\
C_{\dff \bullet}^\iinf\dff(\trf N \trf) 
& 
t_{\dff \bullet}\dff(\trf N\fff,\pff e \trf)
\arrow[r, "\dis \off\qff \tau_{\dff e}"]
\arrow[l, "\dis \off\qff \lambda_{\dff e}"']
&
e_{\bullet}^\lf\dff(\trf X\fff,\pff \mathcal{U}\trf)
\pff.
\end{tikzcd}
\]

\vspace{-7.5pt}
is\dss commutative and\dss leads\sss
to\sss the following\sss commutative diagram of\dss homology\dss
groups\vspace{4.5pt}
\[
\quad
\begin{tikzcd}[column sep=boommm, row sep=boommm]\dis
H_{\dff *}^\iinf\dff(\trf N \trf) 
\arrow[d, "\dis \qff ="]
& 
H_{\dff *}\dff(\trf N\fff,\pff C \trf)
\arrow[d, "\dis \qff \Phi_{\dff *}"]
\arrow[r, red, "\dis \off\qff \tau_{\trf C\dff *}\off", line width=0.8pt]
\arrow[l, "\dis \off\qff \lambda_{\trf C\dff *}"']
&
H_{\dff *}^\lf\dff(\trf X\fff,\pff \mathcal{U}\trf)
\arrow[d, "\dis \qff \varphi_{\dff *}"]
\\
H_{\dff *}^\iinf\dff(\trf N \trf) 
& 
H_{\dff *}\dff(\trf N\fff,\pff e \trf)
\arrow[r, "\dis \off\qff \tau_{\dff e\dff *}\off"]
\arrow[l, red, "\dis \off\qff \lambda_{\dff e\dff *}"', line width=0.8pt]
&
\widetilde{H}_{\dff *}^\lf\dff(\trf X\fff,\pff \mathcal{U}\trf)
\pff,
\end{tikzcd}
\]

\vspace{-7.5pt}
where\dss
$H_{\dff *}\dff(\trf N\fff,\pff e \trf)$\dss
denotes\sss the cohomology\sss of\trs the\sss total\sss complex\dss
$t_{\dff \bullet}\dff(\trf N\fff,\pff e \trf)$\nnsp.\oss

The red arrows are isomorphisms.\oss
Indeed,\oss
since\dss the covering\dss 
$\mathcal{U}$\dss is\dss $e_{\bullet}$\dnsp-acyclic,\pss
$\lambda_{\dff e\dff *}$\sss
is\dss an\dss isomorphism\dss by\trs Lemma\qss \ref{acyclic-coverings},\oss
and\trs Lemma\qss \ref{homology-columns-are-exact}\qss implies\sss that\sss
$\tau_{\qff C\dff *}$\sss is\dss always an\sss isomorphism,\oss
as we already\dss pointed out.\oss
By\dss inverting\dss these\sss two arrows 
we get\dss the commutative diagram\vspace{4.5pt}
\[
\quad
\begin{tikzcd}[column sep=boommm, row sep=boommm]\dis
H_{\dff *}^\iinf\dff(\trf N \trf) 
\arrow[d, "\dis \qff ="]
& 
H_{\dff *}\dff(\trf N\fff,\pff C \trf)
\arrow[d, "\dis \qff \Phi_{\dff *}"]
\arrow[l]
&
H_{\dff *}^\lf\dff(\trf X\fff,\pff \mathcal{U}\trf)
\arrow[l, red, line width=0.8pt]
\arrow[d, "\dis \qff \varphi_{\dff *}"]
\\
H_{\dff *}^\iinf\dff(\trf N \trf)
\arrow[r, red, line width=0.8pt] 
& 
H_{\dff *}\dff(\trf N\fff,\pff e \trf)
\arrow[r]
&
\widetilde{H}_{\dff *}^\lf\dff(\trf X\fff,\pff \mathcal{U}\trf)
\pff.
\end{tikzcd}
\]

\vspace{-7.5pt}
It\dss follows\dss that\dss 
$H_{\dff *}^\lf\dff(\trf X\fff,\pff \mathcal{U}\trf)
\qff \ttoo\qff
\widetilde{H}_{*}^\lf\dff(\trf X\fff,\pff \mathcal{U}\trf)$\dss
factors\sss through\dss the
canonical\dss homomorphism\vspace{4.5pt}
\[
\quad
\begin{tikzcd}[column sep=large, row sep=boom]\dis
H_{\dff *}^\lf\dff(\trf X\fff,\pff \mathcal{U}\trf) 
\arrow[r]
& 
H_{\dff *}\dff(\trf N\fff,\pff C \trf)
\arrow[r]
&
H_{\dff *}^\iinf\dff(\trf N \trf)\pff.
\end{tikzcd}
\]

\vspace{-7.5pt}
The\sss theorem\dss follows.\oss  \eproof

\newpage
\mysection{Compactly\qss finite\qss and\qss $l_{\dff 1}$\dnsp-homology}{cf-homology}

\myuppar{Singular $l_{\dff 1}$\dnsp-chains.}
Recall\dss that\dss for a\sss topological\sss space\sss $Y$\sss
we denote\sss by\sss
$L_{\trf q}\fff(\trf Y\trf)$\sss the real\sss vector space
of\trs infinite singular $q$\dnsp-chains in\sss $Y$\sss
having\sss finite $l_{\dff 1}$\dnsp-norm.\oss
These chains are called\dss \emph{$l_{\dff 1}$\dnsp-chains}\pss of\dss dimension $q$\nnsp.\oss
There are obvious inclusions\dss
$C_{\dff q}\dff(\trf Y\trf)\qff \subset\pff L_{\trf q}\fff(\trf Y\trf)$\nnsp.\oss
The boundary\sss maps in\sss $C_{\dff \bullet}\dff(\trf Y\trf)$\sss
extend\sss by\sss continuity\dss to
$L_{\trf \bullet}\dff(\trf Y\trf)$\nnsp,\pss
turning\sss
$L_{\trf \bullet}\dff(\trf Y\trf)$
into a complex.\oss
Its homology\sss are denoted\sss by\sss
$H_{\dff *}^{\lone}\fff(\trf Y\trf)$\nnsp.\oss
In\dss this section\sss we will\sss apply\dss the\sss theory\sss of\trs
Section\qss \ref{homological-leray-infinite}\qss to\sss
$e_{\dff \bullet}\off =\off L_{\trf \bullet}$\nsp.\oss

Let\sss $\mathcal{U}$ 
be a star\dss finite covering\sss of\pss $X$\nnsp.\oss
The complex\sss $L_{\trf \bullet}^\lf\trf(\trf X\fff,\pff \mathcal{U}\trf)$\sss
admits a description similar\sss to\sss the definition of\dss
$C_{\dff \bullet}^\lf\trf(\trf X\fff,\pff \mathcal{U}\trf)$\nnsp.\oss
Namely,\oss suppose\sss that\dss for every\sss $U\qff \in\qff \mathcal{U}$\sss
an\sss $l_{\dff 1}$\dnsp-chain\sss
$\gamma_{\dff U}\qff \in\qff L_{\trf q}\fff(\trf U\trf)$\sss
is\dss given.\oss
Then,\oss since\sss $\mathcal{U}$\dss is\dss
a star\sss finite,\oss the sum\vspace{1.5pt}
\[
\quad
\gamma
\off =\off
\sum\nolimits_{\qff U\qff \in\qff \mathcal{U}}\qff
\gamma_{\dff U}
\]

\vspace{-9pt}
is\dss a\sss well\sss defined\dss infinite chain.\oss
Let\sss
$\mathcal{L}_{\dff q}^{\dff \lf}\trf(\trf X\fff,\pff \mathcal{U}\trf)$\sss
be\sss the vector space of\dss such chains.\oss

The inclusions\sss
$L_{\trf q}\fff(\trf U\trf)
\qff \ttoo\qff
C_{\dff q}^\iinf\dff(\trf X\trf)$\sss
lead\dss to a map\vspace{3pt}\vspace{0.625pt}
\[
\quad
\delta_{\trf 0}\dff \colon\dff
C_{\trf 0}^\iinf\dff(\trf N\fff,\pff L_{\dff q} \trf)
\qff \ttoo\qff
C_{\dff q}^\iinf\dff(\trf X\trf)
\]

\vspace{-9pt}\vspace{0.625pt}
having\dss
$\mathcal{L}_{\dff q}^\lf\trf(\trf X\fff,\pff \mathcal{U}\trf)$\sss
as\sss the image.\oss
Let\vspace{3pt}\vspace{0.625pt}\vspace{-0.125pt}
\[
\quad
\overline{\delta}_{\trf 0}\dff \colon\dff
C_{\trf 0}^\iinf\dff(\trf N\fff,\pff L_{\dff q} \trf)
\qff \ttoo\qff
\mathcal{L}_{\dff q}^{\dff \lf}\trf(\trf X\fff,\pff \mathcal{U}\trf)
\pff
\]

\vspace{-9pt}\vspace{0.625pt}
be\sss the map resulting\sss from changing\sss the\sss target\sss of\dss
$\delta_{\trf 0}$\nsp.\oss

\mypar{Lemma.}{homology-columns-are-exact-lone}
\emph{The\qss following\dss sequence\dss is\dss exact\dff:}\vspace{-3pt}
\[
\quad
\begin{tikzcd}[column sep=large, row sep=normal]\dis
0
&
\mathcal{L}_{\dff q}^{\dff \lf}\trf(\trf X\fff,\pff \mathcal{U}\trf) 
\arrow[l]
&
C_{\trf 0}^\iinf\dff(\trf N\fff,\pff L_{\dff q} \trf) 
\arrow[l, "\dis \off \overline{\delta}_{\trf 0}\qff"']
& 
C_{\trf 1}^\iinf\dff(\trf N\fff,\pff L_{\dff q} \trf)
\arrow[l, "\dis \off \delta_{\trf 1}\qff"'] 
&   
\off \ldots\off . 
\arrow[l, "\dis \off \delta_{\trf 2}"'] 
\end{tikzcd}
\]

\vspace{-9pt}
\proof
The proof\dss is\dss completely\sss similar\sss to\sss the proof\dss
of\qss Lemma\qss \ref{homology-columns-are-exact}.\oss
On\sss the chains of\trs the form\sss $s\dff *\dff \sigma$\sss
the chain\sss homotopy\sss $k_{\trf p}$\sss is\dss defined as before.\oss
The fact\dss that\sss $\mathcal{U}$\sss is\dss star\sss finite ensures\sss
that\dss this definition extends\sss to\sss $l_{\dff 1}$\dnsp-chains.\oss
The homotopy\dss identity\dss holds on\sss 
the chains of\trs the form\sss $s\dff *\dff \sigma$\sss
by\dss the same reason as before
and\dss hence holds on all\sss $l_{\dff 1}$\dnsp-chains.\oss  \eproof

\mypar{Corollary.}{cokernel}
\emph{The map\qss $\overline{\delta}_{\trf 0}$\sss induces an\dss
isomorphism\pss 
$L_{\trf \bullet}^\lf\trf(\trf X\fff,\pff \mathcal{U}\trf) 
\qff \ttoo\qff
\mathcal{L}_{\dff q}^\lf\trf(\trf X\fff,\pff \mathcal{U}\trf)$\nnsp.}\qss \eproof\dnsp

\mypar{Lemma.}{cf-chains}
\emph{If\pss $\mathcal{U}$ is\dss a star\sss finite proper covering,\oss
then a small\sss chain\dss is\dss $\mathcal{U}$\dnsp-finite\dss
if\trs and\dss only\trs if\qss it\dss is\dss compactly\dss finite.}

\proof
Let\dss us\sss prove\sss the\qss ``if''\qss part\dss first.\oss
Given a small\sss chain $\gamma$\nnsp,\oss 
let\dss us\sss write\dss it\sss as a formal\sss sum\vspace{1.5pt}
\[
\quad
\gamma
\off =\off 
\sum\nolimits_{\dff i\dff \in\dff I}\qff a_{\dff i}\dff \sigma_{\dff i}
\]

\vspace{-10.5pt}
of\dss small\sss simplices $\sigma_{\dff i}$\sss with coefficients\sss
$a_{\dff i}\off \neq\off 0$\nnsp.\oss
We may\sss assume\sss that\sss 
$\sigma_{\dff i}\off \neq\off \sigma_{j}$\sss
if\trs $i\off \neq\off j$\nnsp.\oss 
For every\sss $i\qff \in\pff I$\sss
let\sss us\sss choose\sss $U\dff(\dff i\trf)\qff \in\qff \mathcal{U}$\sss 
such\dss that\sss $\sigma_{\dff i}$\sss is\dss
a simplex\sss in\sss $U\dff(\dff i\trf)$\nnsp.\oss
For\dss $U\qff \in\qff \mathcal{U}$\sss let\vspace{1.5pt}
\[
\quad
\gamma_{\dff U}
\off =\off 
\sum\nolimits_{\qff U\dff(\dff i\trf)\qff =\pff U}\qff a_{\dff i}\dff \sigma_{\dff i}
\pff.
\]

\vspace{-10.5pt}
Clearly,\pss
$\gamma$ is\dss equal\dss to\sss the sum of\dss chains\sss $\gamma_{\dff U}$\nsp.\oss
If\dss $\gamma$\sss is\dss a compactly\sss finite chain,\oss
then every\sss $\gamma_{\dff U}$\sss is\dss a\sss finite chain,\oss 
and\dss hence\sss $\gamma$\sss is\dss $\mathcal{U}$\dnsp-finite.\oss

Let\dss us\sss prove\sss the\qss ``only\dss if''\qss part.\oss
If\trs $\gamma$\sss is\dss a $\mathcal{U}$\dnsp-finite chain,\oss
then\sss
$\gamma
\off =\off
\sum\nolimits_{\qff U\qff \in\qff \mathcal{U}}\qff
\gamma_{\dff U}$\nnsp,\oss
where each\sss $\gamma_{\dff U}$\sss is\dss a\sss finite chain\sss in\sss $U$\nnsp.\oss
If\qss
$Z\qff \subset\qff X$\sss is\dss compact,\oss
then\sss $Z$\sss is\dss contained\sss in\dss the union of\dss finitely\sss
many\sss sets\sss $U\qff \in\qff \mathcal{U}$\nnsp.\oss
Since\sss $\mathcal{U}$\sss is\dss star\sss finite,\oss
this implies\sss that\dss $Z$\sss intersects only\sss finitely\sss many\sss sets\sss
$U\qff \in\qff \mathcal{U}$\nnsp.\oss
Since each\sss $\gamma_{\dff U}$\sss is\dss a\sss finite chain,\oss
it\dss follows\sss that\sss only\sss finitely\sss many\sss simplices entering\sss $\gamma$\sss
with non-zero coefficients intersect\sss $Z$\nnsp.\oss
Hence\sss $\gamma$\sss is\dss compactly\sss finite.\oss  \eproof\vspace{0.125pt}

\mypar{Lemma.}{cf-eilenberg}
\emph{If\pss $\mathcal{U}$ is\dss a star\sss finite proper covering,\oss
then\sss the map\qss
$H_{\dff *}^\lf\trf(\trf X\fff,\pff \mathcal{U}\trf)
\qff \ttoo\qff
H_{\dff *}^\cf\trf(\trf X\trf)$\sss
induced\dss by\dss the inclusion\qss
$C_{\dff \bullet}^\lf\trf(\trf X\fff,\pff \mathcal{U}\trf)
\qff \ttoo\qff
C_{\dff \bullet}^\cf\trf(\trf X\trf)$\sss 
is\dss an\dss isomorphism.\oss}\vspace{0.125pt}

\proof
By\qss Lemma\qss \ref{cf-chains}\qss
the complex\sss
$C_{\dff \bullet}^\lf\trf(\trf X\fff,\pff \mathcal{U}\trf)$\sss
is\dss equal\dss to\sss the subcomplex of\dss small\sss chains
of\trs the complex\sss $C_{\dff \bullet}^\cf\trf(\trf X\trf)$\nnsp.\oss
Hence\sss the\sss lemma\dss is\dss an analogue for compactly\sss finite chains of\dss
a classical\dss theorem of\qss Eilenberg\qss \cite{e}\qss about\dss finite chains.\oss
Eilenberg's\dss proof\trs is\dss presented\sss in\sss many\dss textbooks\qss
(see,\oss for example,\oss \cite{sp},\oss Theorem\qss 4.4.14),\oss
although\dss is\dss rarely\sss attributed\dss to\dss Eilenberg.\oss
Eilenberg's\trs proof\dss applies\sss to our situation\sss without\sss any changes.\oss  \eproof

\myuppar{Comparing\dss compactly\dss finite\sss and $l_{\dff 1}$\dnsp-chains.}
As\dss in\trs Section\qss \ref{homological-leray-infinite},\oss
let\qss 
$\widetilde{H}_{\dff *}^\lf\dff(\qff X\fff,\pff \mathcal{U} \trf)$\sss
be\sss the homology of\trs the complex
{\nsp}$L_{\dff \bullet}^\lf\dff(\qff X\fff,\pff \mathcal{U} \trf)$\nnsp.\oss
There\dss is\dss a canonical\dss comparison\sss homomorphism\vspace{3pt}
\begin{equation}
\label{comparison-lone-lf}
\quad
\varphi_{\dff *}\dff \colon\dff
H_{\dff *}^\lf\trf(\trf X\fff,\pff \mathcal{U}\trf)
\qff \ttoo\qff
\widetilde{H}_{\dff *}^\lf\dff(\qff X\fff,\pff \mathcal{U} \trf)
\pff.
\end{equation}

\vspace{-9pt}
If\trs the assumptions of\qss Lemmas\qss \ref{cf-chains}\qss
and\qss \ref{cf-eilenberg}\qss hold,\oss
we can\sss interpret\sss $\varphi_{\dff *}$ as a homomorphism\vspace{3pt}
\begin{equation}
\label{comparison-lone}
\quad
\varphi_{\dff *}\dff \colon\dff
H_{\dff *}^\cf\trf(\trf X\trf)
\qff \ttoo\qff
\widetilde{H}_{\dff *}^\lf\dff(\qff X\fff,\pff \mathcal{U} \trf)
\pff
\end{equation}

\vspace{-9pt}
and\dss the canonical\dss homomorphisms\qss (\ref{leray-map})\qss
as\sss a\dss homomorphism\dss
$l_{\dff \mathcal{U}}\dff \colon\dff
H_{\dff *}^\cf\trf(\trf X\trf)
\qff \ttoo\qff
H_{\dff *}^\iinf\dff(\trf N \trf)$\nnsp.\oss

\mypar{Lemma.}{lone-factorization}
\emph{Suppose\sss that\qss $\mathcal{U}$ is\dss a star\sss finite proper\sss covering.\oss
If\qss $\mathcal{U}$ is\qss $l_{\dff 1}$\dnsp-acyclic,\oss
then\dss the comparison homomorphisms\pss \textup{(\ref{comparison-lone-lf})}\qss
and\oss \textup{(\ref{comparison-lone})}\qss
can\sss be\trs factored\dss through\dss $l_{\dff \mathcal{U}}$\nsp.\oss}

\proof
Since $L_{\dff \bullet}$\dnsp-acyclicity\dss is\dss the same as
$l_{\dff 1}$\dnsp-acyclicity,\oss
this\sss follows\sss 
from\trs Theorem\qss \ref{e-acyclic-open}.\oss  \eproof

\mypar{Theorem.}{norm-zero-acyclic}
\emph{Suppose\sss that\qss $\mathcal{U}$ is\dss a star\sss finite proper covering,\oss
and\dss that\qss $\mathcal{U}$ is\dss countable and\qss $l_{\dff 1}$\dnsp-acyclic.\oss
If\qss a\sss homology\sss class\dss
$h\qff \in\pff H_{\dff *}^\cf\trf(\trf X\trf)$\sss
belongs\sss to\sss the kernel\dss of\qss $l_{\dff \mathcal{U}}$\nsp,\oss
then\dss $\norm{h}\off =\off 0$\nnsp.\oss}

\proof
Suppose\sss that\dss
$h\qff \in\pff H_{\dff n}^\cf\trf(\trf X\trf)$\dss
belongs\sss to\sss the kernel\dss of\qss $l_{\dff \mathcal{U}}$\nsp.\oss
By\qss Lemma\qss \ref{cf-eilenberg}\qss the homology\sss class $h$
can\sss be represented\dss by\sss a cycle\dss
$\gamma\qff \in\qff C_{\dff n}^\lf\trf(\trf X\fff,\pff \mathcal{U}\trf)$\nnsp.\oss
Lemma\qss \ref{lone-factorization}\qss implies\sss
that\dss\vspace{3pt}
\[
\quad
\varphi_{\dff *}\dff(\dff h\trf)
\off =\off
0
\off \in\off
\widetilde{H}_{\dff n}^\lf\dff(\qff X\fff,\pff \mathcal{U} \trf)
\pff.
\]

\vspace{-9pt}
By\trs Corollary\qss \ref{cokernel}\fff\qss the map\dss $\overline{\delta}_{\trf 0}$\sss
induces an\dss isomorphism\sss between
$\widetilde{H}_{\dff \bullet}^\lf\dff(\qff X\fff,\pff \mathcal{U} \trf)$
and\dss the homology\sss of\trs the complex
$\mathcal{L}_{\dff \bullet}^{\dff \lf}\trf(\trf X\fff,\pff \mathcal{U}\trf)$\nnsp.\oss
It\dss follows\sss that\dss the inclusion\dss
$C_{\dff n}^\lf\trf(\trf X\fff,\pff \mathcal{U}\trf)
\qff \ttoo\qff
\mathcal{L}_{\dff n}^{\dff \lf}\trf(\trf X\fff,\pff \mathcal{U}\trf)$\dss
takes\sss the cycle $\gamma$ representing\sss $h$\sss to a boundary.\oss
In\sss other\dss terms,\qss\vspace{3pt}
\[
\quad
\gamma
\off =\off 
\partial\trf \beta
\]

\vspace{-9pt}
for\sss some\trs 
$\beta\qff \in\qff
\mathcal{L}_{\dff n\dff +\dff 1}^{\dff \lf}\trf(\trf X\fff,\pff \mathcal{U}\trf)$\nnsp.\oss
By\dss the\sss definition\sss of\trs
$\mathcal{L}_{\dff n\dff +\dff 1}^{\dff \lf}\trf
(\trf X\fff,\pff \mathcal{U}\trf)$\nnsp,\oss\vspace{3pt}\vspace{0.75pt}
\[
\quad
\beta
\off =\off
\sum\nolimits_{\qff U\qff \in\qff \mathcal{U}}\pff \beta_{\dff U}
\pff
\]

\vspace{-9pt}\vspace{0.75pt}
for some chains\dss
$\beta_{\dff U}\qff \in\qff L_{\dff n\dff +\dff 1}\fff(\trf U\trf)$\nnsp.\oss
Let\dss us\sss choose an arbitrary\dss $\varepsilon\qff >\qff 0$\nnsp.\oss
Since\sss $\mathcal{U}$\sss is\dss countable,\oss
there exists a family\sss of\dss real\dss numbers\dss
$\varepsilon_{\dff U}\qff >\qff 0$\nnsp,\qss
$U\qff \in\qff \mathcal{U}$\nnsp,\dff\oss such\dss that\qss\vspace{3pt}\vspace{0.75pt}
\[
\quad
\sum\nolimits_{\qff U\qff \in\qff \mathcal{U}}\qff \varepsilon_{\dff U}
\off =\off 
\varepsilon
\pff.
\]

\vspace{-9pt}\vspace{0.75pt}
Since\dss $\norm{\fff\beta_{\dff U}}\qff <\qff \infty$\dss
for every\sss $U$\nnsp,\oss
every\sss chain\sss $\beta_{\dff U}$\sss can\sss be presented as a sum\dss
$\beta_{\dff U}\off =\off \alpha_{\dff U}\qff +\qff \omega_{\dff U}$\dss
of\trs two chains\sss 
$\alpha_{\trf U}\dff,\off \omega_{\dff U}
\qff \in\qff L_{\dff n\dff +\dff 1}\fff(\trf U\trf)$\sss
such\dss that\sss
$\alpha_{\trf U}$\sss is\dss finite and\dss
$\norm{\fff\omega_{\dff U}}
\qff <\qff 
\varepsilon_{\dff U}$\nnsp.\oss
Let\vspace{4.5pt}
\[
\quad
\alpha
\off =\off
\sum\nolimits_{\qff U\qff \in\qff \mathcal{U}}\qff \alpha_{\trf U}
\quad\
\mbox{and}\quad\
\omega
\off =\off
\sum\nolimits_{\qff U\qff \in\qff \mathcal{U}}\qff \omega_{\dff U}
\pff.
\]

\vspace{-7.5pt}
Then\sss
$\alpha
\qff \in\pff 
C_{\dff n\dff +\dff 1}^\lf\dff(\trf X\fff,\pff \mathcal{U}\trf)
\qff \subset\qff
C_{\dff n\dff +\dff 1}^\cf\dff(\trf X\trf)$\sss
and\dss
$\omega
\qff \in\qff 
\mathcal{L}_{\dff n\dff +\dff 1}^{\dff \lf}\trf
(\trf X\fff,\pff \mathcal{U}\trf)$\nnsp.\oss
Now,\pss
$\gamma
\off =\off 
\partial\trf \beta$\dss
implies\sss that\vspace{3pt}
\[
\quad
\gamma
\off =\off
\partial\dff \alpha\qff +\qff \partial\dff \omega
\]

\vspace{-9pt}
and\dss hence\dss 
$\gamma\qff -\qff \partial\dff \alpha
\off =\off
\partial\dff \omega$\nnsp.\oss
Therefore\vspace{4.5pt}
\[
\quad
\norm{\gamma\qff -\qff \partial\dff \alpha}
\off =\off
\norm{\fff\partial\dff \omega}
\qff \leq\qff 
(\dff n\qff +\qff 1\dff)\dff \norm{\omega}
\]

\vspace{-33pt}
\[
\quad
\phantom{\norm{\gamma\qff -\qff \partial\dff \alpha}
\off =\off
\norm{\fff\partial\dff \omega}
\qff }
\leq\qff
(\dff n\qff +\qff 1\dff)\dff
\sum\nolimits_{\qff U\qff \in\qff \mathcal{U}}\qff \norm{\omega_{\dff U}}
\]

\vspace{-33pt}
\[
\quad
\phantom{\norm{\gamma\qff -\qff \partial\dff \alpha}
\off =\off
\norm{\fff\partial\dff \omega}
\qff }
<\qff
(\dff n\qff +\qff 1\dff)\dff
\sum\nolimits_{\qff U\qff \in\qff \mathcal{U}}\pff \varepsilon_{\dff U}
\off =\off 
(\dff n\qff +\qff 1\dff)\dff
\varepsilon
\pff.
\]

\vspace{-7.5pt}
Hence $h$\sss can be represented\trs by\sss chains 
with arbitrarily\sss small\dss norm,\oss
i.e.\qss $\norm{h}\off =\off 0$\nnsp.\oss  \eproof

\mypar{Theorem.}{norm-finite-almost-acyclic}
\emph{Suppose\sss that\qss $\mathcal{U}$ is\dss a star\sss finite proper covering,\oss
and\dss that\qss $\mathcal{U}$ is\dss countable and\sss almost\qss $l_{\dff 1}$\dnsp-acyclic.\oss
If\pss 
$h\qff \in\pff H_{\dff *}^\cf\trf(\trf X\trf)$\sss
belongs\sss to\sss the kernel\dss of\qss $l_{\dff \mathcal{U}}$\nsp,\oss
then\dss $\norm{h}\off <\off \infty$\nnsp.\oss}

\proof
Let\dss $U_{\fff e}\qff \in\qff \mathcal{U}$\sss be\sss the exceptional\sss set,\oss
the one which\dss is\dss allowed\sss not\dss to be\sss $l_{\dff 1}$\dnsp-acyclic.\oss
By\qss Lemma\qss \ref{cf-eilenberg}\qss the homology\sss class $h$
can\sss be considered as an element\sss of\dss
$H_{\dff n}^\lf\dff(\trf X\fff,\pff \mathcal{U}\trf)$\sss
and\dss represented\dss by a chain\dss
$\gamma
\qff \in\qff 
C_{\dff n}^\lf\dff(\trf X\fff,\pff \mathcal{U}\trf)$\sss
for some $n$\nnsp.\oss
Let\dss us consider\sss the commutative diagram\vspace{-1.5pt}
\[
\quad
\begin{tikzcd}[column sep=boommmm, row sep=boommmm]\dis
H_{\dff \bullet}^\iinf\dff(\trf N \trf) 
\arrow[d, "\dis \qff ="]
& 
H_{\dff *}\dff(\trf N\fff,\pff C \trf)
\arrow[d, "\dis \qff \Phi_{\dff *}"]
\arrow[l, "\dis \off\qff \lambda_{\trf C\dff *}"']
&
H_{\dff *}^\lf\dff(\trf X\fff,\pff \mathcal{U}\trf)
\arrow[l, red, "\dis \off\qff \tau_{\trf C\dff *}^{\dff -\dff 1}"', line width=0.8pt]
\arrow[d, "\dis \qff \varphi_{\dff *}"]
\\
H_{\dff *}^\iinf\dff(\trf N \trf) 
& 
H_{\dff *}\dff(\trf N\fff,\pff L \trf)
\arrow[r, "\dis \off\qff \tau_{\dff L\dff *}\off"]
\arrow[l, red, "\dis \off\qff \lambda_{\dff L\dff *}"', line width=0.8pt]
&
\widetilde{H}_{\dff *}^\lf\dff(\trf X\fff,\pff \mathcal{U}\trf)
\end{tikzcd}
\]

\vspace{-7.5pt}
similar\sss to\sss the diagrams used\sss in\sss the proof\dss of\trs
Theorem\qss \ref{e-acyclic-open}.\oss
Since\sss $\mathcal{U}$\sss is\dss not\sss assumed\dss to be
$l_{\dff 1}$\dnsp-acyclic,\oss the homomorphism\sss 
$\lambda_{\dff L\dff *}$\sss may\sss be not\sss invertible.\oss
But\sss since\sss
$l_{\dff \mathcal{U}}\dff(\dff h\trf)\off =\off 0$\nnsp,\oss\vspace{3pt}
\[
\quad
\widetilde{h}
\off =\off
\Phi_{\dff *}\dff\left(\trf \tau_{\trf C\dff *}^{\dff -\dff 1}\dff(\dff h\trf)\qff\right)
\off \in\off
H_{\dff n}\dff(\trf N\fff,\pff L \trf)
\]

\vspace{-9pt}
belongs\sss to\sss the kernel\sss of\dss $\lambda_{\dff L\dff *}$\nsp.\oss
Hence\trs Lemma\qss \ref{double-complex}\qss
implies\sss that\sss $\widetilde{h}$\sss belongs to\sss the image of\vspace{3pt}
\[
\quad
H_{\dff n}\dff(\qff C_{\trf 0}^\iinf\dff(\dff N\fff,\pff L_{\dff \bullet} \dff)\trf)
\qff \ttoo\qff
H_{\dff n}\dff(\trf N\fff,\pff L \trf)
\pff.
\]

\vspace{-9pt}
Since among\sss sets\sss in\sss $\mathcal{U}$\sss only\sss $U_{\fff e}$\sss
can\sss be not\sss $l_{\dff 1}$\dnsp-acyclic,\pss
$\widetilde{h}$\trs belongs\sss to\sss the image of\trs the homology of\trs
the summand\sss $L_{\dff \bullet}\dff(\trf U_{\fff e}\trf)$\sss of\trs
$C_{\trf 0}^\iinf\dff(\dff N\fff,\pff L_{\dff \bullet} \dff)$\nnsp.\oss
It\dss follows\dss that\dss there exists an 
$l_{\dff 1}$\dnsp-cycle\vspace{3pt}
\[
\quad
\gamma\fff'
\qff \in\qff
L_{\dff n}\dff(\trf U_{\fff e}\trf)
\off \subset\off
C_{\trf 0}^\iinf\dff(\dff N\fff,\pff L_{\dff n} \dff)
\]

\vspace{-9pt}
such\sss that\dss $\widetilde{h}$\trs is\dss equal\dss to\sss the $l_{\dff 1}$\dnsp-homology\sss
class of\dss $\gamma\fff'$\sss and\dss hence\sss
$\tau_{\dff L\dff *}\dff(\trf \widetilde{h}\qff)$\sss
is\dss equal\dss to\sss the homology\sss class of\trs the cycle\sss $\gamma'$\sss
considered as an element\sss of\trs
$\mathcal{L}_{\dff \bullet}^\lf\dff(\qff X\fff,\pff \mathcal{U} \trf)$\nnsp.\oss

Now\sss the commutativity\sss of\trs the right\sss square of\trs the
above diagram\dss implies\sss that\sss
$\varphi_{\dff *}\dff(\dff h\trf)$\sss
is\dss equal\dss to\sss the homology\sss class $h\fff'$\sss of\dss $\gamma\fff'$\nnsp.\oss 
It\dss follows\dss that\qss
$\gamma\qff -\qff \gamma'$\qss
is\dss a\sss boundary\dss in\dss
$\mathcal{L}_{\dff \bullet}^\lf\dff(\qff X\fff,\pff \mathcal{U} \trf)$\nnsp.\oss
As\sss in\dss the proof\dss of\qss Theorem\qss \ref{norm-zero-acyclic},\oss
for every\dss $\varepsilon\qff >\qff 0$\dss
there exist\sss chains\dss
$\alpha
\qff \in\pff 
C_{\dff n\dff +\dff 1}^\lf\dff(\trf X\fff,\pff \mathcal{U}\trf)$\sss
and\dss
$\omega
\qff \in\qff 
\mathcal{L}_{\dff n\dff +\dff 1}^{\dff \lf}\trf
(\trf X\fff,\pff \mathcal{U}\trf)$\sss
such\dss that\sss $\norm{\omega}\qff <\qff \varepsilon$\sss
and\vspace{3pt}
\[
\quad
\gamma\qff -\qff \gamma\fff'
\off =\off
\partial\dff \alpha\qff +\qff \partial\dff \omega
\pff,
\]

\vspace{-9pt}
i.e.\qss 
$\gamma\qff -\qff \partial\dff \alpha
\off =\off
\gamma\fff'\qff +\qff \partial\dff \omega$\nnsp.\oss
The cycle\sss $\gamma\qff -\qff \partial\dff \alpha$\sss
is\dss $\mathcal{U}$\dnsp-finite\sss and\dss represents\sss
$h$\nnsp.\oss
On\dss the other\sss hand\sss $\norm{\gamma\fff'}\qff <\qff \infty$\sss
and\sss $\norm{\omega}\qff <\qff \varepsilon$\nnsp,\oss
and\dss hence\sss
$\norm{\gamma\qff -\qff \partial\dff \alpha}\qff <\qff \infty$\nnsp.\oss
Therefore\sss
$\norm{h}\qff <\qff \infty$\nnsp.\oss  \eproof

\mysection{Extensions\qss of\pss coverings\qss and\qss $l_{\dff 1}$\dnsp-homology}{extensions}

\myuppar{Extensions of\dss coverings.}
Let\sss $\mathcal{U}$\sss be a covering of\dss $X$\sss and\sss $X\fff'$\sss 
be a\sss space containing\sss $X$\nnsp.\oss
Recall\qss (see\qss \cite{i3},\oss Section\qss 4)\qss 
that\sss an\qss \emph{extension of\dss $\mathcal{U}$\sss to\sss $X\fff'$}\trs
is\dss a map\sss $U\qff \ttoo\qff U\fff'$\sss assigning\sss to every\sss
$U\qff \in\qff \mathcal{U}$\sss a subset\sss $U\fff'\qff \subset\qff X\fff'$\sss
such\dss that\sss $U\fff'\dff \cap\dff X\off =\off U$\sss in\sss such a way\sss that\dss
the collection\sss $\mathcal{U}\fff'$ of\trs the sets $U\fff'$\sss 
is\dss a covering of\dss $X\fff'$\nnsp.\oss
The set\sss $\mathcal{U}\fff'$\sss uniquely\sss determines\sss the extension,\oss
i.e.\qss the map\sss $U\qff \ttoo\qff U\fff'$\sss and\sss is\dss identified\sss with\sss it.\oss
There\dss is\dss an obvious simplicial\dss map\sss
$N_{\trf \mathcal{U}}\qff \ttoo\qff N_{\trf \mathcal{U}\fff'}$\nnsp.\oss

The extension\sss $\mathcal{U}\fff'$\sss is\dss said\dss to be\qss \emph{nerve-preserving}\pss
if\trs this map\dss is\dss a simplicial\dss isomorphism,\oss
which\dss is\dss then\sss treated as\sss the identity.\oss
If\dss $\sigma$\sss is\dss a simplex of\dss $N_{\trf \mathcal{U}}$\nsp,\oss
then\sss $\num{\sigma}\fff'$\sss denotes\sss the intersection of\dss
sets\sss $U\fff'$\sss corresponding\dss to\sss the vertices of\dss $\sigma$\nnsp.\oss
Clearly,\pss 
$\num{\sigma}\qff \ttoo\qff \num{\sigma}\fff'$\qss
is\dss a map\sss
$\mathcal{U}^{\dff \cap}\qff \ttoo\qff \mathcal{U}\fff'^{\dff \cap}$\dnsp.\oss

Suppose\sss that\sss $\mathcal{U}$\sss is\dss a weakly $l_{\dff 1}$\dnsp-acyclic
covering of\dss $X$\nnsp.\oss
By\qss \cite{i3},\oss Corollary\qss 4.2,\oss
there exists a space\dss $X\fff'\qff \supset\qff X$\sss
and a nerve-preserving\sss extension\sss $\mathcal{U}\fff'$\sss of\trs
$\mathcal{U}$\sss to\sss $X\fff'$\sss such\dss that\sss
$\mathcal{U}\fff'$\sss is\sss $l_{\dff 1}$\dnsp-acyclic
and\dss the inclusion\dss $X\qff \ttoo\qff X'$\sss induces an\sss isomorphism
of\trs the fundamental\dss groups.\oss
Therefore\dss $X\qff \ttoo\qff X'$\sss
induces isometric isomorphisms\sss in\sss bounded cohomology.\oss
By\sss a\sss theorem of\qss Cl.\dss L\"{o}h\qss \cite{l}\qss 
this implies\sss that\dss $X\qff \ttoo\qff X'$\sss
induces isometric isomorphisms\sss in\sss $l_{\dff 1}$\dnsp-homology.\oss

Moreover\halfff,\oss if\dss $\mathcal{U}$\sss is\dss open,\oss 
then\sss $\mathcal{U}\fff'$\sss
can\sss be assumed\dss to be also open.\oss
See\qss \cite{i3},\oss Corollary\qss 4.2.\oss
The same argument\sss shows\sss that\dss if\trs the interiors\sss
$\inte U$\nnsp,\qss $U\qff \in\qff \mathcal{U}$\nnsp,\pss
cover\sss $X$\nnsp,\oss 
then\sss one can assume\sss that\dss 
the interiors\sss $\inte U\fff'$\sss cover\sss $X\fff'$\nnsp.\oss
Also,\oss one can assume\sss that\sss $\inte U\fff'\qff \supset\qff \inte U$\nnsp,\oss
where\sss the first\sss interior\dss is\dss taken\sss in\sss $X\fff'$\sss
and\dss the second\sss in\sss $X$\nnsp.\oss
See\qss \cite{i3},\oss the end of\trs the proof\dss of\qss Theorem\qss 4.1.\oss
We will\dss need also\sss the following\sss simple property\sss of\trs
the construction of\dss $\mathcal{U}\fff'$\dnsp.\oss

\mypar{Lemma.}{retraction}
\emph{Let\qss $U\qff \longmapsto\qff U_{\dff +}$
be\sss the map\sss 
establishing\sss that\dss the covering\sss
$\mathcal{U}$\sss is\dss weakly\sss $l_{\dff 1}$\dnsp-acyclic.\oss
Then\dss there exists a retraction\sss
$r\dff \colon\dff
X\fff'\qff \ttoo\qff X$\sss
such\dss that\sss
$r\trf(\trf U\fff'\trf)\qff \subset\pff U^{\dff +}$\sss 
for every\sss $U\qff \in\qff \mathcal{U}^{\dff \cap}$\dnsp.\oss}

\proof
For\sss every\sss $U\qff \in\qff \mathcal{U}^{\dff \cap}$\sss
the subset\sss $U\fff'\qff \subset\qff X\fff'$\sss
is\dss obtained\dss from $U$ by\sss attaching discs
along\sss loops contractible\sss in $X$ and\dss then attaching\sss some\qss
``collars''\qss to ensure\sss that\dss $\inte U\fff'\qff \supset\qff \inte U$\nnsp.\oss
See\qss \cite{i3},\oss the proof\dss of\qss Theorem\qss 4.1.\oss
Under\sss our assumptions\sss the\sss loops used\dss to attach discs\sss to $U$
are contractible in\sss $U^{\dff +}$\dnsp.\oss
It\dss follows\sss that\dss there exists a retraction\sss
$r\dff \colon\dff
X\fff'\qff \ttoo\qff X$\sss
such\dss that\sss
$r\trf(\trf U\fff'\trf)\qff \subset\pff U^{\dff +}$\sss 
for every\sss $U\qff \in\qff \mathcal{U}$\dnsp.\oss  \eproof

\mypar{Theorem.}{norm-zero-weakly-acyclic}
\emph{Suppose\sss that\dss $\mathcal{U}$ is\dss a star\sss finite proper covering\sss
and\dss that\sss $\mathcal{U}$ is\dss countable and\dss weakly $l_{\dff 1}$\dnsp-acyclic.\oss
If\pss
$h\qff \in\pff H_{\dff *}^\cf\trf(\trf X\trf)$\sss
belongs\sss to\sss the kernel\dss of\qss $l_{\dff \mathcal{U}}$\nsp,\oss
then\dss $\norm{h}\off =\off 0$\nnsp.}

\proof
Let\sss $\mathcal{U}\fff'$\sss be\sss the extension of\trs the covering\sss
$\mathcal{U}$\sss to a space\sss $X\fff'\qff \supset\qff X$\sss as above.\oss
The closures of\trs the sets\sss $U\fff'$\nnsp, where\sss $U\qff \in\qff \mathcal{U}$\nnsp,\oss
may\sss be not\sss compact\dss because,\oss in\dss general,\pss
$U\fff'$\sss is\dss obtained\sss from\sss $U$\sss by\sss attaching an\sss
infinite number of\dss discs and\qss ``collars''.\oss
Therefore\sss the covering\sss $\mathcal{U}\fff'$\sss is\dss not\dss proper\sss
in\sss general,\oss and we cannot\sss apply\trs
Theorem\qss \ref{norm-zero-acyclic}\qss to\sss it.\oss

But\dss the proof\dss of\qss Theorem\qss \ref{norm-zero-acyclic}\qss
applies with only\sss minor\sss modifications.\oss
Let\dss us\sss represent\dss the homology\sss class $h$\sss by\sss a cycle\vspace{3pt}
\[
\quad
\gamma\qff \in\qff C_{\dff n}^\lf\trf(\trf X\fff,\pff \mathcal{U}\trf)
\off \subset\off
C_{\dff \bullet}^\lf\trf(\trf X\fff'\fff,\pff \mathcal{U}\fff'\trf)
\pff.
\]

\vspace{-9pt}
The cycle\sss
$\gamma$\sss defines also a\sss homology\sss class\sss
$h\fff'\qff \in\pff
H_{\dff n}^\lf\trf(\trf X\fff'\fff,\pff \mathcal{U}\fff'\trf)$\nnsp.\oss
Clearly,\oss the diagram\vspace{0pt}
\[
\quad
\begin{tikzcd}[column sep=boomm, row sep=boomm]\dis
H_{\dff *}^\lf\dff(\trf X\fff,\pff \mathcal{U}\trf)
\arrow[r, "\dis l_{\dff \mathcal{U}}"] 
\arrow[d]
& 
H_{\dff *}^\iinf\dff(\trf N \trf)
\arrow[d, "\dis \dff ="]
\\
H_{\dff *}^\lf\dff(\trf X\fff'\fff,\pff \mathcal{U}\fff'\trf)
\arrow[r, "\dis l_{\dff \mathcal{U}\fff'}"] 
& 
H_{\dff *}^\iinf\dff(\trf N \trf)
\end{tikzcd}
\]

\vspace{-12pt}
is\dss commutative.\oss
Since\sss 
$l_{\dff \mathcal{U}}\dff(\dff h\trf)\off =\off 0$\nnsp,\oss
this\sss implies\dss that\sss
$l_{\dff \mathcal{U}\fff'}\dff(\dff h\fff'\trf)\off =\off 0$\nnsp.\oss
Together\dss with\trs Theorem\qss \ref{e-acyclic-open}\qss this\sss implies\sss
that\sss $h\fff'\off =\off 0$\nnsp,\oss
i.e.\qss
$\gamma$\sss is\dss a boundary\dss in\sss the chain complex\sss
$L_{\trf \bullet}^\lf\trf(\trf X\fff'\fff,\pff \mathcal{U}\fff'\trf)$\nnsp.\oss

In\sss view\sss of\qss Corollary\qss \ref{cokernel}\qss this\sss implies\sss
that\sss $\gamma$\sss is\dss a boundary\sss in\dss 
$\mathcal{L}_{\trf \bullet}^\lf\trf(\trf X\fff'\fff,\pff \mathcal{U}\fff'\trf)$\nnsp,\oss
i.e.\qss\vspace{1.5pt}\vspace{-0.29pt}
\[
\quad
\gamma
\off =\off 
\partial\trf \beta
\pff
\]

\vspace{-10.5pt}\vspace{-0.29pt}
for\dss some\trs 
$\beta\qff \in\qff
\mathcal{L}_{\dff n\dff +\dff 1}^{\dff \lf}\trf(\trf X\fff'\fff,\pff \mathcal{U}\fff'\trf)$\nnsp.\oss
Then\vspace{3pt}\vspace{-0.58pt}
\[
\quad
\beta
\off =\off
\sum\nolimits_{\qff U\qff \in\qff \mathcal{U}}\pff \beta_{\dff U}
\pff
\]

\vspace{-9pt}
for some chains\dss
$\beta_{\dff U}\qff \in\qff L_{\dff n\dff +\dff 1}\fff(\trf U\fff'\trf)$\nnsp.\oss
Arguing\sss as\sss in\sss the proof\dss of\qss Theorem\qss \ref{norm-zero-acyclic},\oss
let\dss us\sss choose an arbitrary\dss $\varepsilon\qff >\qff 0$\sss and\sss
represent\sss $\varepsilon$\sss as a sum\sss
$\varepsilon
\off =\off
\sum\nolimits_{\qff U\qff \in\qff \mathcal{U}}\qff 
\varepsilon_{\dff U}$\sss 
of\dss positive numbers\sss $\varepsilon_{\dff U}\qff >\qff 0$\nnsp.\oss
Next,\oss let\dss us\sss represent\sss each\sss $\beta_{\dff U}$\sss as a sum\dss
$\beta_{\dff U}\off =\off \alpha_{\dff U}\qff +\qff \omega_{\dff U}$\dss
of\dss chains\sss 
$\alpha_{\dff U}\dff,\off \omega_{\dff U}
\qff \in\qff L_{\dff n\dff +\dff 1}\fff(\trf U\trf)$\sss
such\dss that\sss
$\alpha_{\dff U}$\sss is\dss finite and\dss
$\norm{\fff\omega_{\dff U}}
\qff <\qff 
\varepsilon_{\dff U}$\nnsp.\oss
Let\vspace{3pt}\vspace{-0.58pt}
\[
\quad
\omega
\off =\off
\sum\nolimits_{\qff U\qff \in\qff \mathcal{U}}\qff \omega_{\dff U}
\quad\
\mbox{and}\quad\
\alpha
\off =\off
\sum\nolimits_{\qff U\qff \in\qff \mathcal{U}}\qff \alpha_{\dff U}
\pff.
\]

\vspace{-9pt}
For\sss every\sss $U\qff \in\qff \mathcal{U}$\sss
the chain\sss $r_{\dff *}\dff(\trf \alpha_{\dff U}\trf)$\sss
is\dss a chain\sss in\sss $U^{\dff +}$\dnsp.\oss
Since\sss the family\sss of\dss sets\sss $U^{\dff +}$\sss 
is\dss compactly\sss finite,\oss
the infinite chain\vspace{3pt}
\[
\quad
r_{\dff *}\dff(\trf \alpha\trf)
\off =\off
\sum\nolimits_{\qff U\qff \in\qff \mathcal{U}}\qff
r_{\dff *}\dff(\trf \alpha_{\dff U}\trf)
\]

\vspace{-9pt}
is\dss well\sss defined and compactly\sss finite.\oss
Since $r$\sss is\dss a retraction,\oss\vspace{3pt}
\[
\quad
\gamma
\qff -\qff
\partial\trf\bigl(\trf
r_{\dff *}\dff(\trf \alpha\trf)
\trf\bigr)
\off =\off
r_{\dff *}\dff(\trf \gamma\trf)
\qff -\qff
\partial\trf\bigl(\trf
r_{\dff *}\dff(\trf \alpha\trf)
\trf\bigr)
\]

\vspace{-36pt}
\[
\quad
\phantom{\gamma
\qff -\qff
\partial\trf\bigl(\trf
r_{\dff *}\dff(\trf \alpha\trf)
\trf\bigr)
\off }
=\off
r_{\dff *}\dff(\trf \gamma\trf)
\qff -\qff
r_{\dff *}\dff\bigl(\trf \partial\dff\alpha\trf\bigr)
\off =\off
r_{\dff *}\dff\bigl(\trf \gamma
\qff -\qff
\partial\dff\alpha \trf\bigr)
\pff.
\]

\vspace{-9pt}\vspace{-0.75pt}
But\sss
$\gamma
\qff -\qff
\partial\dff\alpha
\off =\off
\partial\dff\omega$\sss
and\sss
$\norm{\partial\dff\omega}
\qff \leq\qff
(\dff n\qff +\qff 1\dff)\dff \norm{\omega}
\qff <\qff 
(\dff n\qff +\qff 1\dff)\dff \varepsilon$\nnsp.\oss
It\dss follows\dss that\vspace{4.5pt}
\[
\quad
\norm{r_{\dff *}\dff(\trf \gamma
\qff -\qff
\partial\dff\alpha \trf)}
\off \leq\off
\norm{\gamma
\qff -\qff
\partial\dff\alpha}
\off <\off
(\dff n\qff +\qff 1\dff)\dff \varepsilon
\pff.
\]

\vspace{-7.5pt}
Hence\sss $h$\sss can\sss be represented\trs by\sss chains 
with arbitrarily\sss small\dss norm,\oss
i.e.\qss $\norm{h}\off =\off 0$\nnsp.\oss  \eproof\vspace{0.2pt}

\mypar{Theorem.}{norm-finite-almost-weakly-acyclic}
\emph{Suppose\sss that\dss $\mathcal{U}$ is\dss a star\sss finite proper covering\sss
and\dss that\sss $\mathcal{U}$ is\dss countable and\dss almost\dss weakly $l_{\dff 1}$\dnsp-acyclic.\oss
If\pss
$h\qff \in\pff H_{\dff *}^\cf\trf(\trf X\trf)$\sss
belongs\sss to\sss the kernel\dss of\qss $l_{\dff \mathcal{U}}$\nsp,\oss
then\dss $\norm{h}\off <\off \infty$\nnsp.}\vspace{0.2pt}

\proof
The proof\dss differs from\dss the proof\dss
of\qss Theorem\qss \ref{norm-zero-weakly-acyclic}\qss
is\dss the same way\sss as\sss the proof\dss
of\qss Theorem\qss \ref{norm-finite-almost-acyclic}\qss
differs from\sss the proof\dss of\qss Theorem\qss \ref{norm-zero-acyclic}.\oss
We\sss leave\sss the details\sss to\sss the reader\halfff.\oss  \eproof

\myuppar{Compactly\sss amenable families.}
The\dss
$l_{\dff 1}$\dnsp-acyclicity\sss is\dss implied\dss by\sss a stronger\sss property,\oss 
namely,\oss the amenability.\oss
Suppose\sss that\sss $Z\qff \subset\qff Y\qff \subset\qff X$\sss
and\sss $Z$\sss is\dss path connected.\oss
The subset\sss $Z$\sss is\dss said\dss to be\qss
\emph{amenable\sss in}\dss $Y$\dss
if\trs the image of\trs the map\dss
$\pi_{\dff 1}\dff(\trf Z\fff,\qff z \trf)
\trf \ttoo\trf 
\pi_{\dff 1}\dff(\trf Y\fff,\qff z \trf)$\trs
is\sss amenable.\oss

A family\sss $\mathcal{U}$\sss of\dss subsets of\dss $X$\sss is\dss said\dss to be\qss
\emph{compactly\sss amenable}\pss 
if\dss for every\sss $U\qff \in\qff \mathcal{U}$\sss
a subset\sss $U_{\dff +}\qff \subset\qff X$\sss is\dss given,\oss
such\dss that\sss $U\qff \subset\qff U_{\dff +}$\nsp,\oss
the set\sss $U$\sss is\dss amenable\sss in\sss $U_{\dff +}$\nsp,\oss
and\dss the family\sss of\dss subsets\sss $U_{\dff +}$\sss
is\dss compactly\sss finite.\oss
A covering\sss $\mathcal{U}$\sss is\dss said\dss to be\qss 
\emph{compactly\sss amenable}\pss if\trs it\dss is\dss
compactly\sss amenable as a\sss family\sss and elements of\dss $\mathcal{U}^{\dff \cap}$\sss
are path connected.\oss

Since subgroups of\dss amenable groups are amenable,\oss
if\dss $Z\fff'$\sss is\dss a\sss path connected subset\sss 
of\dss a set\sss $Z$\sss amenable\sss in\sss $Y$\dnsp,\oss
then\sss $Z\fff'$\sss is\dss also amenable in\sss $Y$\nnsp.\oss
At\dss the same\sss time amenable groups are $l_{\dff 1}$\dnsp-acyclic.\oss
It\dss follows\sss that\sss a compactly\sss amenable covering\dss
is\dss compactly\sss $l_{\dff 1}$\dnsp-acyclic.\oss\vspace{0.2pt}

\mypar{Theorem.}{norm-zero-amenable}
\emph{Suppose\sss that\dss $\mathcal{U}$ is\dss a star\sss finite proper covering\sss
and\sss $\mathcal{U}$ is\dss countable and\dss compactly\sss amenable.\oss
If\pss
$h\qff \in\pff H_{\dff *}^\cf\trf(\trf X\trf)$\sss
belongs\sss to\sss the kernel\dss of\qss $l_{\dff \mathcal{U}}$\nsp,\oss
then\dss $\norm{h}\off =\off 0$\nnsp.}\vspace{0.2pt}

\proof
In view of\trs the remarks preceding\dss the\sss theorem,\oss
this follows\sss from\trs Theorem\qss \ref{norm-zero-weakly-acyclic}.\oss  \eproof

\myuppar{Almost\dss compactly\dss amenable\sss families.}
A\sss family\sss $\mathcal{U}$\sss of\dss subsets of\dss $X$\sss 
is\dss said\sss to be\qss
\emph{almost\dss compactly\sss amenable}\pss 
if\dss subsets\sss $U_{\dff +}\qff \subset\qff X$\sss
with\dss the same properties\sss as in\sss the definition
of\dss compactly\sss amenable families 
are given\sss for\sss every subset\sss
$U\qff \in\qff \mathcal{U}$\dnsp,\oss
except\halfff,\pss perhaps,\oss of\dss a single exceptional\sss set\sss 
$U_{\fff e}\qff \in\qff \mathcal{U}$\dnsp.\oss
A covering\sss $\mathcal{U}$\sss is\dss said\dss to be\qss 
\emph{almost\sss compactly\sss amenable}\pss if\trs it\dss is\dss
almost\sss compactly\sss amenable as a\sss family\sss 
and elements of\dss $\mathcal{U}^{\dff \cap}$\sss
are path connected,\oss except,\oss perhaps,\oss the set\sss $U_{\fff e}$\nsp.\oss
Clearly,\oss an almost\sss compactly\sss amenable covering\dss is\dss
almost\dss weakly $l_{\dff 1}$\dnsp-acyclic.\oss\vspace{0.2pt}

\mypar{Theorem.}{norm-finite-almost-compactly-amenable}
\emph{Suppose\sss that\dss $\mathcal{U}$ is\dss a star\sss finite proper covering\sss 
and\sss $\mathcal{U}$ is\dss countable and\dss almost\sss compactly\sss amenable.\oss
If\pss
$h\qff \in\pff H_{\dff *}^\cf\trf(\trf X\trf)$\sss
belongs\sss to\sss the kernel\dss of\qss $l_{\dff \mathcal{U}}$\nsp,\oss
then\dss $\norm{h}\off <\off \infty$\nnsp.}\vspace{0.2pt}

\proof
In view of\trs the remarks preceding\dss the\sss theorem,\oss
this follows\sss from\trs Theorem\qss \ref{norm-finite-almost-weakly-acyclic}.\oss  \eproof

\myuppar{Families compactly\sss amenable\sss in\dss the sense of\qss Gromov.}
Gromov's\trs definitions
of\dss amenable subsets and\sss coverings are slightly\sss different.\oss
Let\sss us\sss say\dss that\sss $Z$\sss is\qss
\emph{amenable\sss in\sss the sense of\qss Gromov}\pss in\sss $Y$\sss
if\dss every\sss path connected component\sss of\dss $Z$\sss 
is\dss amenable in\sss $Y$\sss in our sense.\oss

A family\sss $\mathcal{U}$\sss is\qss
\emph{compactly\sss amenable\sss in\sss the sense of\qss Gromov}\pss
if\dss for every\sss $U\qff \in\qff \mathcal{U}$\sss
a subset\sss $U_{\dff +}\qff \subset\qff X$\sss is\dss given,\oss
such\dss that\sss $U\qff \subset\qff U_{\dff +}$\nsp,\oss
the set\sss $U$\sss is\dss amenable in\sss the sense of\trs Gromov\sss 
in\sss $U_{\dff +}$\nsp,\oss
and\dss the family\sss of\dss subsets\sss $U_{\dff +}$\sss
is\dss compactly\sss finite.\oss

A\sss family\sss $\mathcal{U}$\sss 
is\qss
\emph{almost\dss compactly\sss amenable\sss in\sss the sense of\qss Gromov}\pss 
if\dss subsets\sss $U_{\dff +}\qff \subset\qff X$\sss
with\dss the same properties\sss as in\sss the definition
of\dss compactly\sss amenable families\sss in\sss the sense of\qss Gromov 
are given\sss for\sss every\sss
$U\qff \in\qff \mathcal{U}$\dnsp,\oss
except\halfff,\pss perhaps,\oss of\dss finitely\dss many exceptional\sss sets.\oss
Gromov\qss \cite{gro}\qss uses\sss the\sss term\qss
\emph{sequence\qss ``amenable''\trs at\dss infinity}\pss
for\sss a\sss slightly\sss different\dss notion.\oss

The main difference of\trs Gromov's\dss version of\trs these notions\dss
is\dss the\sss lack of\dss any\sss assumptions of\dss path connectedness.\oss
Still,\oss a\sss large part\sss of\dss our\sss theory\sss
survives in\sss this context.\oss

Let\dss us\sss relax\dss the assumption of\dss 
$e_{\dff \bullet}$\dnsp-acyclicity\sss of\trs the covering\sss $\mathcal{U}$\sss
in\dss Section\qss \ref{homological-leray-infinite}\qss
by\sss the assumption\sss that\dss the homology\sss groups of\dss complexes\sss
$e_{\dff \bullet}\dff(\trf Z\trf)$\dss with\sss
$Z\qff \in\qff \mathcal{U}^{\dff \cap}$\sss 
vanish\sss in dimensions\dss $>\qff 0$\nnsp.\oss
Then,\oss in order\sss to keep\dss Lemma\qss \ref{acyclic-coverings},\oss
we need\dss to replace\sss the complex\sss $C_{\dff \bullet}^\iinf\dff(\trf N \trf)$\sss
by\dss the cokernel\dss of\trs 
the horizontal\dss boundary\sss operator\sss\vspace{3pt}\vspace{-0.75pt}
\[
\quad
d_{\dff 1}\dff \colon\dff
c_{\trf \bullet}\fff(\trf N\fff,\pff e_{\trf 1} \trf)
\qff \ttoo\qff
c_{\trf \bullet}\fff(\trf N\fff,\pff e_{\trf 0} \trf)
\pff.
\]

\vspace{-9pt}\vspace{-0.75pt}
We will\sss denote\sss this cokernel\dss by\sss $C_{\dff \bullet}$\nsp,\oss
and\sss denote by\sss $H_{\dff *}$\sss its homology\sss groups.\oss
These homology\sss groups play\sss now\sss the role of\dss $H_{\dff *}^\iinf\dff(\trf N \trf)$\nnsp.\oss
With\sss these changes\sss the arguments of\qss Section\qss \ref{homological-leray-infinite}\qss
still\dss work and\sss show\sss 
that\dss the comparison\sss homomorphism\qss\vspace{3pt}\vspace{-0.75pt}
\[
\quad
\varphi_{\dff *}\dff \colon\dff
H_{\dff p}^\lf\trf(\trf X\fff,\pff \mathcal{U}\trf)
\off \ttoo\off
\widetilde{H}_{\dff p}^\lf\dff(\qff X\fff,\pff \mathcal{U} \trf)
\]

\vspace{-9pt}\vspace{-0.75pt}
factors\sss through a canonical\dss map\dss
$H_{\dff *}^\lf\trf(\trf X\fff,\pff \mathcal{U}\trf)
\off \ttoo\off
H_{\dff *}$\nsp.\oss
Clearly,\pss 
$C_{\dff p} \off =\off 0$\sss
if\dss $p\qff >\qff \dim\trf N$\nnsp.\oss
Therefore\sss
$H_{\fff p} \off =\off 0$\sss
for\dss $p\qff >\qff \dim\trf N$\nnsp.\oss
It\dss follows\dss that\dss
$\varphi_{\dff *}$\sss is\dss equal\dss to $0$\sss 
for\dss $p\qff >\qff \dim\trf N$\nnsp.\oss

Now\sss for a homology\sss class\sss
$h\qff \in\pff H_{\dff p}^\cf\trf(\trf X\trf)$\sss
we can\sss require\sss that\sss
$p\qff >\qff \dim\trf N$\sss
instead of\dss requiring\dss that\sss $h$\sss
belongs\sss to\sss the kernel\dss of\qss $l_{\dff \mathcal{U}}$\nsp.\oss

\mypar{Theorem.}{norm-zero-amenable-gromov}
\emph{Let\dss $\mathcal{U}$ is\dss a star\sss finite proper covering\sss 
which\dss is\dss countable and\dss 
compactly\sss amenable in\sss the sense of\pss Gromov.\oss
If\pss
$h\qff \in\pff H_{\dff p}^\cf\trf(\trf X\trf)$\sss
and\dss $p\qff >\qff \dim\trf N$\nnsp,\oss
then\dss $\norm{h}\off =\off 0$\nnsp.\oss}  \eproof

\mypar{Theorem.}{norm-finite-almost-compactly-amenable-gromov}
\emph{Let\dss $\mathcal{U}$ is\dss a star\sss finite proper covering\sss 
which\dss is\dss countable and\dss 
almost\dss compactly\sss amenable in\sss the sense of\pss Gromov.\oss
If\pss
$h\qff \in\pff H_{\dff p}^\cf\trf(\trf X\trf)$\sss
and\dss $p\qff >\qff \dim\trf N$\nnsp,\oss
then\dss $\norm{h}\off <\off \infty$\nnsp.}

\proof
In\dss this situation only\sss one additional\sss step\dss is\dss needed.\oss
Namely,\oss one needs\sss to replace exceptional\sss sets by\dss their 
union.\oss  \eproof

\newpage
\mysection{Removing\pss weakly\dss $l_{\dff 1}$\dnsp-acyclic\qss subspaces}{removing}

\myuppar{The restriction\sss homomorphisms.}
Let\sss $X$\sss be a\sss topological\sss space,\oss
and\dss let\sss $Y\qff \subset\qff X$\sss be a closed subset.\oss
Let\dss us\sss construct\sss some chain\dss maps\vspace{1.5pt}
\[
\quad
r_{\dff \smallsetminus\dff Y}
\qff \colon\dff
C_{\dff \bullet}^\cf\dff(\trf X\trf)
\qff \ttoo\qff
C_{\dff \bullet}^\cf\dff(\trf X\qff \smallsetminus\qff Y\trf)
\pff.
\]

\vspace{-12pt}
Let\dss $\sigma\dff \colon\dff \Delta^n\qff \ttoo\qff X$\sss
be a singular $n$\dnsp-simplex.\oss
If\dss $\sigma\dff(\dff \Delta^n\trf)\qff \subset\qff Y$\nnsp,\oss
then\sss $\ry\trf(\dff \sigma\dff)\off =\off 0$\nnsp.\oss
Otherwise $\sigma^{\dff -\dff 1}\dff(\trf  X\qff \smallsetminus\qff Y\trf)$\sss
is\dss a non-empty\sss open subset\sss of\dss $\Delta^n$\dnsp.\oss
Let\dss us\sss triangulate\sss this subset\sss by\sss some geometric\qss
(rectilinear)\qss simplices and\sss linearly\sss order\sss the vertices of\trs
this\sss triangulation.\oss
Then every\sss $n$\dnsp-dimensional\sss simplex $\alpha$ of\trs the\sss triangulation
defines an affine singular simplex\sss 
$\alpha'$\sss in\sss $\Delta^n$\sss
and a singular $n$\dnsp-simplex\sss
$\sigma\dff \circ\trf \alpha'
\dff \colon\dff
\Delta^n\qff \ttoo\qff
X\qff \smallsetminus\qff Y$\dnsp.\oss
Let\sss $\ry\trf(\dff \sigma\dff)$\sss be\sss the sum of\dss all\dss these
singular simplices\sss $\sigma\dff \circ\trf \alpha'$\nnsp.\oss
The chain\sss $\ry\trf(\trf \sigma\dff)$\sss is\sss compactly\sss finite
because a compact\sss subset\sss of\dss 
$\sigma^{\dff -\dff 1}\dff(\trf  X\qff \smallsetminus\qff Y\trf)$\sss
intersects only\sss a finite number of\dss simplices of\trs the\sss triangulation.\oss
The map\sss $\ry$\sss extends\sss by\sss linearity\sss to\sss 
$C_{\dff \bullet}^\cf\dff(\trf X\trf)$\sss and\dss maps compactly\sss finite chains
in\sss $X$\sss to compactly\sss finite chains in\sss $X\qff \smallsetminus\qff Y$\dnsp.\oss

In\sss general,\oss such a map\sss $\ry$\sss does not\sss commute with\dss the
boundary\sss operators,\oss i.e.\qss is\dss not\sss a chain\sss map.\oss
In order\sss to ensure\sss that\sss $\ry$\sss is\dss a chain\sss map
one needs to construct\dss the above\sss triangulations recursively,\oss
starting\sss with\dss the\sss tautological\dss triangulations for $0$\dnsp-simplices.\oss
If\trs triangulations are already\sss constructed\sss for singular\sss $m$\dnsp-simplices
with\sss $m\qff <\qff n$\sss and $\sigma$\sss is\dss a singular $n$\dnsp-simplex,\oss
then\sss 
$\partial\dff \Delta^n
\dff \cap\trf
\sigma^{\dff -\dff 1}\dff(\trf  X\qff \smallsetminus\qff Y\trf)$\sss
is\dss already\sss triangulated and\sss one can extend\sss this\sss
triangulation\sss to\sss a\sss triangulation of\dss
$\sigma^{\dff -\dff 1}\dff(\trf  X\qff \smallsetminus\qff Y\trf)$\nnsp.\oss
By\sss continuing\sss in\sss this way\sss we will\sss get\sss
a map\sss $\ry$\sss commuting\sss with\sss the boundary\sss operators.\oss

The resulting\dss map $\ry$ depends on\sss the choice 
of\trs these\sss triangulations.\oss
But\sss different\sss choices\sss led\dss to chain-homotopic chain\sss maps.\oss
Given\sss two choices of\trs triangulations,\oss a chain\sss homotopy\dss
between\dss the corresponding\sss maps\sss $\ry$\sss can\sss be constructed
similarly\sss to\sss the maps\sss $\ry$\sss themselves.\oss
Namely,\oss for every singular $n$\dnsp-simplex $\sigma$\sss
the\sss two\sss triangulations of\dss 
$\sigma^{\dff -\dff 1}\dff(\trf  X\qff \smallsetminus\qff Y\trf)$\sss
can\sss be considered as a\sss triangulation of\dss
$\sigma^{\dff -\dff 1}\dff(\trf  X\qff \smallsetminus\qff Y\trf)
\dff \times\dff
\{\trf 0\fff,\qff 1\trf\}$\nnsp,\oss
and one needs\sss to extend\dss this\sss triangulation\sss to\sss
$\sigma^{\dff -\dff 1}\dff(\trf  X\qff \smallsetminus\qff Y\trf)
\dff \times\dff
[\trf 0\fff,\qff 1\trf]$\nnsp.\oss
If\trs these\sss extensions are constructed\sss recursively,\oss
then\sss they\sss will\sss define a chain\sss homotopy\sss 
between\sss two maps\sss $\ry$\nsp.\oss
Therefore\sss the map\vspace{1.5pt}
\[
\quad
\hry\dff \colon\dff
H_{\dff \bullet}^\cf\dff(\trf X\trf)
\qff \ttoo\qff
H_{\dff \bullet}^\cf\dff(\trf X\qff \smallsetminus\qff Y\trf)
\pff.
\]

\vspace{-12pt}
induced\sss by\sss $\ry$\sss does not\sss depend on\dss the choice of\trs triangulations.\oss
We will\sss assume\sss that\sss a choice of\trs triangulations\dss
is\dss fixed and $\ry$\sss is\dss 
the corresponding\dss map.\oss\vspace{-0.125pt}

\mypar{Lemma.}{zero}
\emph{Let\dss $A$\sss be a\dss topological\sss space and\dss
$n\qff \in\pff \nnn$\nnsp.\oss
Suppose\sss that\dss
$H^{\dff l_{\dff 1}}_{\dff n}\dff(\trf A\trf)
\off =\off 0$\nnsp.\oss
Then\dss there exists a constant\dss $K$\dss 
with\dss the following\sss property.\oss
If\dss $z$\dss is\dss an $l_{\dff 1}$\dnsp-cycle
in\sss $A$\nnsp,\oss
then\sss $z\off =\off \partial\dff u$\sss
for some $l_{\dff 1}$\dnsp-chain $u$ such\dss that\dss
$\norm{u}\qff \leq\qff K\dff \norm{z}$\nnsp.\oss}

\proof
This\sss observation\dss is\dss due\sss to\dss Matsumoto\dss and\trs Morita\qss \cite{mm}.\oss
Let\dss us\sss consider\dss the boundary\sss operator\sss
$\partial\dff \colon\dff
C^{\dff l_{\dff 1}}_{\dff n\dff +\dff 1}\dff(\trf A\trf)
\qff \ttoo\qff
C^{\dff l_{\dff 1}}_{\dff n}\dff(\trf A\trf)$\nnsp.\oss
By\dss the assumption,\oss its\sss image\dss is\dss the subspace\sss
$Z^{\dff l_{\dff 1}}_{\dff n}\dff(\trf A\trf)$\sss of\dss cycles.\oss
Since\sss the\sss latter\dss 
is\dss defined as\sss the kernel\sss of\dss a bounded operator,\oss
it\dss is\dss closed and\dss hence\dss is\dss a\dss Banach\dss space.\oss
The kernel\sss of\trs $\partial$\sss is\dss the subspace of\dss cycles
$Z^{\dff l_{\dff 1}}_{\dff n\dff +\dff 1}\dff(\trf A\trf)$\nnsp.\oss
Therefore\sss $\partial$\sss induces a\sss linear\sss isomorphism\vspace{1.5pt}\vspace{-0.25pt}
\[
\quad
\partial\fff'\dff \colon\dff
C^{\dff l_{\dff 1}}_{\dff n\dff +\dff 1}\dff(\trf A\trf)\left/\dff
Z^{\dff l_{\dff 1}}_{\dff n\dff +\dff 1}\dff(\trf A\trf)\right.
\qff \ttoo\qff
Z^{\dff l_{\dff 1}}_{\dff n}\dff(\trf A\trf)
\pff.
\]

\vspace{-10.5pt}\vspace{-0.25pt}
Since\sss $\partial\fff'$\sss is\dss bounded,\oss
the open\sss mapping\dss theorem\sss implies\sss that\sss its inverse\dss
is\dss bounded.\oss
Let\dss $K\fff'$\dss be\sss the norm of\trs the inverse.\oss
Then\sss for every\sss 
$z\qff \in\qff Z^{\dff l_{\dff 1}}_{\dff n}\dff(\trf A\trf)$\sss
there exists an element\sss $u'$\sss of\dss the above quotient\sss
such\dss that\dss 
$z\off =\off \partial\fff'\fff u'$\dss and\dss
the norm of\sss $u'$\sss is\dss
$\leq\qff K\fff'\dff\norm{z}$\nsp.\oss
By\dss the definition of\trs the norm on a quotient\sss of\dss a\dss
Banach\dss space\sss by\sss a closed subspace,\oss 
for every\sss $\varepsilon\qff >\qff 0$\sss
there\dss is\dss a representative\sss $u$\sss of\dss $u'$\sss
such\dss that\dss
$\norm{u}
\qff \leq\qff
K\fff'\dff\norm{z}
\qff +\qff
\varepsilon$\nnsp.\oss
If\dss $\norm{z}\qff >\qff 0$\nnsp,\oss
we can\sss take\sss 
$\varepsilon
\off =\off 
\norm{z}$\nsp.\oss
If\dss $\norm{z}\off =\off 0$\nnsp,\oss
then\sss $z\off =\off 0$\sss and\sss $u\off =\off 0$\sss
is\dss a representative of\dss $u'$\nnsp.\oss
In\sss both cases\dss
$\norm{u}
\qff \leq\qff
(\trf K\fff'\qff +\qff 1\dff)\dff\norm{z}$\nnsp.\oss
Since\dss
$\partial\fff u
\off =\off
\partial\fff'\dff u'
\off =\off
z$\nnsp,\oss
we can\sss take\sss 
$K\off =\off K\fff'\qff +\qff 1$\nnsp.\oss  \eproof

\mypar{Lemma.}{double-boundary-zero}
\emph{Let\dss $A$\sss be a\dss topological\sss space and\dss
$n\qff \in\pff \nnn$\nnsp.\oss
Suppose\sss that\dss}\vspace{3pt}\vspace{-0.25pt}
\[
\quad
H^{\dff l_{\dff 1}}_{\dff n}\dff(\trf A\trf)
\off =\off
H^{\dff l_{\dff 1}}_{\dff n\dff -\dff 1}\dff(\trf A\trf)
\off =\off 
0
\pff.
\]

\vspace{-9pt}\vspace{-0.25pt}
\emph{Then\sss for every\dss $\varepsilon\qff >\qff 0$\dss 
there exists\sss $\delta\qff >\qff 0$\sss
with\dss the following\dss property.\qff\oss
If\qss $a$\sss is\dss a finite $n$\dnsp-chain\sss in\sss $A$\sss
and\qss 
$\norm{\partial\dff a}\qff <\qff \delta$\nnsp,\oss
then\dss there exist\sss a\dss finite chain\sss $a'$\sss
such\dss that\dss 
$a'\qff -\qff a\off =\off \partial\fff b$\trs
for some finite chain\sss $b$\sss and\qss
$\norm{a'}\qff <\qff \varepsilon$\nnsp.\oss}

\proof
Let\dss $K$\dss be\sss the constant\dss having\dss the property\sss of\qss
Lemma\qss  \ref{zero}\qss with\sss $n\qff -\qff 1$\sss
in\dss the role of\dss $n$\nnsp.\oss
The boundary\sss $\partial\dff a$\sss is\dss a finite cycle and\dss hence
an $l_{\dff 1}$\dnsp-cycle.\oss
Lemma\qss  \ref{zero}\qss implies\sss that\sss 
$\partial\dff a\off =\off \partial\dff d$\sss 
for\sss some $l_{\dff 1}$\dnsp-chain\sss
$d$\sss such\dss that\dss
$\norm{d\fff}
\qff \leq\qff 
K\dff \norm{\partial\dff a}$\nsp.\oss
Clearly,\pss
$d\qff -\qff a$\sss is\dss an $l_{\dff 1}$\dnsp-cycle
of\dss dimension\sss $n$\nnsp.\oss
By\dss the assumption\sss $d\qff -\qff a$\sss is $l_{\dff 1}$\dnsp-homologous\sss
to $0$\nnsp,\oss
i.e.\sss $d\qff -\qff a\off =\off \partial\dff e$\sss
for some $l_{\dff 1}$\dnsp-chain\sss $e$\nnsp.\oss
Let\dss $\delta\qff >\qff 0$\dss be such\dss that\dss
$K\dff \delta\qff <\qff \varepsilon$\nnsp,\oss
and\dss let\dss us\sss choose\sss $\delta\fff'\qff >\qff 0$\sss
such\dss that\vspace{3pt}\vspace{-0.25pt}
\[
\quad
K\dff \delta
\qff +\qff
(\dff n\qff +\qff 2\dff)\dff \delta\fff'
\off <\off
\varepsilon
\pff.
\]

\vspace{-9pt}\vspace{-0.25pt}
Let\dss us\sss represent\sss $e$\sss
as a\sss sum\sss $e\off =\off b\pff +\qff b\fff'$\dnsp,\oss
where\sss $b$\sss is\dss a\dss finite\dss chain\sss and\sss
$\norm{b\fff'}\qff \leq\off \delta\fff'$\dnsp,\oss
and\dss let\sss 
$a'\off =\off a\qff +\qff \partial\fff b$\nnsp.\oss
Then\sss $a'\qff -\qff a\off =\off \partial\fff b$\nnsp,\qss
$a'$\sss is\dss a\sss finite chain,\oss and\sss
$a'\off =\off d\qff -\qff \partial\fff b\fff'$\dnsp.\oss
Hence\vspace{3pt}\vspace{-0.25pt}
\[
\quad
\norm{a'}
\off \leq\off
\norm{d\fff}
\qff +\qff
\norm{\partial\fff b\fff'}
\off \leq\off
K\dff \norm{\partial\dff a}
\qff +\qff
(\dff n\qff +\qff 2\dff)\dff \norm{b\fff'}
\off <\off
\varepsilon
\pff.
\]

\vspace{-9pt}\vspace{-0.25pt}
The\sss lemma\sss follows.\oss  \eproof

\mypar{Lemma.}{double-boundary}
\emph{Let\dss $C$\sss be a subspace of\dss a\sss topological\sss space\sss $D$\sss 
and\dss let\dss $n\qff \in\pff \nnn$\nnsp.\oss
Suppose\sss that\dss
$C$\sss is\dss path connected\sss and\sss 
weakly $l_{\dff 1}$\dnsp-acyclic\sss in\dss $D$\nnsp.\oss
For every\dss $\varepsilon\qff >\qff 0$\dss there exists\sss $\delta\qff >\qff 0$\sss
with\dss the following\dss property.\oss
If\qss $c$\sss is\dss a finite chain\sss in\sss $C$\sss
and\dss $\norm{\partial\dff c}\qff <\qff \delta$\nnsp,\oss
then\dss there exist\sss a\sss finite chain\sss $c'$\sss
in\dss $D$\sss such\dss that\sss 
$c'\qff -\qff c\off =\off \partial\dff b$\sss
for some finite chain\sss $b$\sss and\qss
$\norm{c'}\qff <\qff \varepsilon$\nnsp.\oss}

\proof
Let\sss $C\fff'$\sss be\sss the result\sss of\dss attaching discs\sss to\sss $C$\sss
along\sss a set\sss of\trs loops such\dss that\dss their\sss homotopy\sss classes
generate\sss the kernel\sss of\trs the inclusion\sss homomorphism\vspace{3pt}\vspace{0.125pt}
\[
\quad
\pi_{\dff 1}\dff(\dff C\fff,\qff x\trf)
\qff \ttoo\qff
\pi_{\dff 1}\dff(\dff D\fff,\qff x\trf)
\pff,
\]

\vspace{-9pt}
where\sss $x\qff \in\qff A$\nnsp.\oss
Then\sss $\pi_{\dff 1}\dff(\dff C\fff'\fff,\qff x\trf)$\sss
is\dss amenable and\dss hence\sss the bounded cohomology\sss
$\widehat{H}^{\trf \bullet}\fff(\trf C\fff'\trf)$\sss are zero.\oss
By\sss a well\dss known\dss theorem of\qss Matsumoto\dss and\trs Morita\qss \cite{mm}\qss
(see also\qss \cite{i3}\qss for a\sss proof\dff)\qss
this implies\sss that\sss
$H^{\dff l_{\dff 1}}_{\dff \bullet}\dff(\trf C\fff'\trf)
\off =\off
0$\nnsp.\oss
By\qss Lemma\qss \ref{double-boundary-zero}\qss there exist\sss
a chain $a$\sss in\sss $C\fff'$\sss such\dss that\dss
$a'\qff -\qff a\off =\off \partial\dff z$\trs
for\sss some finite chain $z$\sss in\sss $C\fff'$\sss and\qss
$\norm{a'}\qff <\qff \varepsilon$\nnsp.\oss
It\dss remains\sss to\sss turn\sss $a'\fff,\pff z$\sss into chains\sss in\dss $D$\sss
while keeping\dss its\sss properties.\oss
Let\dss $D\fff'\off =\off C\fff'\dff \cup\dff D$\nnsp.\oss
Since\sss the\sss loops used\dss to attach discs\sss to $C$
are contractible in\sss $D$\nnsp,\oss
there exists a\sss retraction\dss
$r\dff \colon\dff 
D\fff'\qff \ttoo\qff D$\dnsp.\oss
Let\sss $c'\off =\off r_{\dff *}\dff(\dff a\trf)$\nnsp.\oss
Then\dss
$\norm{c'}
\qff \leq\qff
\norm{a\fff}
\qff <\qff 
\varepsilon$\nnsp,\oss
the chain\sss
$b\off =\off r_{\dff *}\dff(\dff z\trf)$\sss 
is\dss finite,\oss
and\vspace{3pt}
\[
\quad
c'\qff -\qff c
\off =\off
r_{\dff *}\dff(\dff a\trf)\qff -\qff r_{\dff *}\dff(\dff c\trf)
\off =\off
r_{\dff *}\dff(\dff a\qff -\qff c\trf)
\off =\off
r_{\dff *}\dff(\trf \partial\dff z\trf)
\off =\off
\partial\dff r_{\dff *}\dff(\dff z\trf)
\pff.
\]

\vspace{-9pt}
The\sss lemma\sss follows.\oss  \eproof

\myuppar{Supports\sss and\sss parts\sss of\dss singular\sss chains.}
For a singular simplex\sss
$\sigma\dff \colon\dff
\Delta^n\trf \ttoo\qff X$\sss
let\dss
$\overline{\sigma}\off =\off \sigma\dff(\trf \Delta^n\trf)$\nnsp.\oss
The\qss \emph{support}\qss $\supp\dff(\dff c\trf)$\sss 
of\trs the singular chain\qss (\ref{singular-chain})\qss 
is\dss defined as\sss the union\vspace{3pt}
\[
\quad
\supp\dff c\dff(\dff c\trf)
\off =\off
\bigcup\nolimits_{\dis\qff a_{\dff \sigma}\qff \neq\qff 0}\pff \overline{\sigma}
\pff.
\]

\vspace{-9pt}
For a subset\trs $Y\qff \subset\qff X$\dss the\dss \emph{$Y$\dnsp-part}\qss 
$c\trf|\dff Y$\sss
of\trs the singular chain\qss (\ref{singular-chain})\qss is\dss defined as\sss the chain\vspace{3pt}
\[
\quad
c\trf|\trf Y
\off =\off
\sum\nolimits_{\qff \overline{\sigma}\qff \cap\qff Y\qff \neq\qff \varnothing}\pff 
a_{\dff \sigma}\dff \sigma
\pff.
\]

\vspace{-9pt}
If\dss $c$\sss is\dss compactly\sss finite,\oss
then\sss $c\trf|\trf Y$\sss is\dss also compactly\sss finite.\oss 
The\qss \emph{intersection}\qss $c\dff \cap\dff Y$\dss is\vspace{3pt}
\[
\quad
c\dff \cap\dff Y
\off =\off
\sum\nolimits_{\qff \overline{\sigma}\qff \subset\qff Y}\pff 
a_{\dff \sigma}\dff \sigma
\pff.
\]

\vspace{-9pt}
Clearly,\oss
$c\off =\off 
c\trf|\trf Y\qff +\qff c\dff \cap\dff (\trf X\ssm Y\trf)$\nnsp.\oss

\myuppar{Surgery\sss of\trs chains.}
Recall\dss that\sss $Y$\sss is\dss a closed subset\sss of\dss $X$\nnsp.\oss
Suppose\sss that\trs $Z$\sss is\dss a compact\sss component\sss of\pss $Y$\sss
and\dss that\trs $C$\sss be a compact\trs Hausdorff\trs neighborhood\dss 
of\pss $Z$\sss disjoint\dss from\sss $Y\ssm Z$\nnsp.\oss
Suppose\sss that\sss $\gamma$\sss is\dss a compactly\sss finite cycle\sss in\sss $X\ssm Y$\nnsp.\oss
We would\dss like\sss remove\sss from\sss $\gamma$\sss some part\sss contained\sss in\sss $C$\sss
and\dss replace\sss it\dss by\sss a finite chain\sss in\sss $X$\sss 
without\sss noticeably\sss increasing\dss the norm.\oss

More precisely,\oss given\sss $\varepsilon\qff >\qff 0$\nnsp,\oss
we would\dss like\sss to find\sss an open set\sss
$U\qff \subset\qff C$\sss and a finite chain\sss $\psi$\sss in\sss $X$\sss
such\dss that\sss
$\partial\trf (\trf \gamma\dff \cap\dff U\trf)
\off =\off
\partial\dff \psi$\sss
and\sss
$\norm{\psi}\qff <\qff \varepsilon$\nnsp.\oss 
Since\dss
$\gamma
\off =\off
\gamma\trf|\trf X\ssm U
\qff +\qff
\gamma\dff \cap\dff U$\dss
is\dss a cycle,\oss the chain\sss
$\gamma\trf|\trf X\ssm U
\qff +\qff
\psi$\sss
is\dss also a cycle and\dss its\sss norm\dss 
$<\qff \norm{\gamma}\qff +\qff \varepsilon$\nnsp.\oss
If\trs the part\sss $\gamma\dff \cap\dff U$\sss removed\sss 
from\sss $\gamma$\sss and\dss the chain\sss $\psi$\sss
depend only\sss on\sss $\gamma\trf |\trf C$\nnsp,\oss
this operation could\dss be performed\sss for several\sss
components\sss $Z$\sss simultaneously.\oss

In our applications all\sss components of\dss $Y$\sss will\dss be  compact\sss 
and\sss $\gamma$\sss will\dss represent\dss the homology\sss class\sss
$\hry\dff(\dff h\trf)$\dss for some\dss $
h\qff \in\qff H_{\dff \bullet}^\cf\dff(\trf X\trf)$\nnsp.\oss
In\dss this case we would\dss like\sss to get\sss a representative of\dss $h$\sss
after\sss performing\dss this operation for 
all\sss components of\dss $Y$\sss simultaneously.\oss

\mypar{Lemma.}{chain-modification}
\emph{Let\pss $\delta\pff >\pff 0$\nnsp.\oss
Then\dss for\sss every\sss compactly\trs finite cycle\dss $\gamma$\dss
in\pss $X\ssm Y$\qss such\dss that\qss
$\norm{\gamma}\qff <\qff \infty$\qss
there\sss exists\sss an open\dss neighborhood\pss 
$U$ of\pss $Z$\sss contained\dss in\trs $C$\sss and\sss such\dss that}\vspace{3pt}
\[
\quad
\norm{(\trf \partial\trf (\trf \gamma\dff \cap\dff U\trf)}
\off <\off\qff
\delta
\pff.
\]

\vspace{-9pt}
\proof
Since\sss
$\norm{\gamma}\qff <\qff \infty$\nnsp,\oss
for every\sss $\varepsilon\qff >\qff 0$\sss one can\sss write $\gamma$
as a sum\sss
$\gamma\off =\off \gamma\fff'\qff +\qff \gamma\fff''$\sss
of\dss a\sss finite chain\sss $\gamma\fff'$\sss and a chain\sss $\gamma\fff''$\sss
such\dss that\dss $\norm{\gamma\fff''}\qff <\qff \varepsilon$\nnsp.\oss
Since\sss $\gamma\fff'$\sss is\dss finite,\oss
the support\sss $\supp\dff \gamma\fff'$\sss is\dss compact\sss
and\dss hence\sss the intersection\sss 
$C\dff \cap\dff \supp\dff \gamma\fff'$\sss is\dss also compact.\oss
Since\sss $C$\sss is\qss Hausdorff\halfff,\oss this\sss intersection\dss is\dss closed\sss
and\dss its complement\sss $U$\sss in\sss $C$\sss is\dss open.\oss
Clearly,\oss every\sss simplex\sss contained\sss in\sss $U$\sss and entering\sss $\gamma$\sss
with a non-zero coefficient\sss enters $\gamma\fff''$ with\dss the same coefficient.\oss
It\dss follows\sss that\sss
$\norm{\gamma\dff \cap\dff U}
\qff \leq\qff
\norm{\gamma\fff''}
\qff <\qff 
\varepsilon$\nnsp.\oss
Let\sss $n$\sss be\sss the dimension of\dss $\gamma$ and\dss let\sss us\sss take\sss
$\varepsilon\off =\off \delta/(\dff n\qff +\qff 1\dff)$\nnsp.\oss
Then\sss
$\norm{(\trf \partial\trf (\trf \gamma\dff \cap\dff U\trf)}
\qff \leq\qff
(\dff n\qff +\qff 1\dff)\dff \norm{\gamma\dff \cap\dff U}
\qff <\qff
\delta$\nnsp.\oss  \eproof

\myuppar{Families of\dss compact\sss subspaces.}
Now\sss we need\dss to impose further\dss restrictions on\sss $Y$\dnsp.\oss
Suppose\sss that\sss $Y$\sss is\dss presented as\sss the union of\dss a\sss family\sss
$\mathcal{Z}$\sss of\dss pair-wise disjoint\sss compact\sss subspaces of\dss $X$\nnsp.\oss
Suppose\sss further\dss that\dss for\sss every\sss
$Z\qff \in\qff \mathcal{Z}$\sss a compact\sss
neighborhood\dss $C_{\trf Z}$\dss of\dss $Z$\sss is\dss given,\oss
and\dss that\dss the neighborhoods\dss $C_{\trf Z}$\sss are pair-wise disjoint.\oss
Suppose\sss that\sss every\sss $C_{\trf Z}$\sss is\trs Hausdorff\trs and\dss path connected.\oss
Let\sss $V_{\trf Z}$\sss be\sss the interior of\dss $C_{\trf Z}$\sss and\dss
$V$\dss be\sss the union of\dss all\sss sets\sss $V_{\trf Z}$\nsp.\oss

\mypar{Lemma.}{homology}
\emph{Let\sss $h\qff \in\qff H_{\dff n}^\cf\dff(\trf X\trf)$\sss
and\dss let\trs $\gamma$\sss be a compactly\trs finite chain\sss in\dss $X\ssm Y$\sss
representing\dss the homology\sss class\sss $\ry\trf(\dff h\trf)$\nnsp.\oss
Let\qss $X\fff'\qff \subset\qff X\ssm Y$\sss be a\sss closed\dss set\sss
containing\oss 
$X\ssm V$\dnsp.\oss
Then\dss there exists\sss a compactly\trs finite chain\sss $s$\sss in\dss $V$\sss
such\dss that\qss 
$\gamma\trf |\qff X\fff'\pff +\pff s$\dss 
is\dss a cycle representing\sss $h$\nnsp.\oss}

\proof
Let\sss $c\qff \in\qff
C_{\dff n}^\cf\dff(\trf X\trf)$\sss
be a cycle representing\sss $h$\nnsp.\oss
Then\vspace{3pt}
\begin{equation}
\label{gamma}
\quad
\gamma 
\off =\off
\ry\trf(\dff c\trf)
\qff +\qff
\partial\dff \beta
\end{equation}

\vspace{-9pt}
for some chain\dss
$\beta\qff \in\qff
C_{\dff n\dff +\dff 1}^\cf\dff(\trf X\qff \smallsetminus\qff Y\trf)$\nnsp.\oss
Let\vspace{3pt}
\[
\quad
\gamma\fff'
\off =\off
\gamma\trf|\qff X\fff'\qff,\quad\
r\fff'
\off =\off 
\ry\trf(\dff c\trf)\trf|\qff X\fff'\qff,\quad\
\mbox{and}\quad\
\beta\fff'
\off =\off
\beta\trf|\qff X\fff'
\pff.
\]

\vspace{-9pt}
If\dss $\tau$\sss is\dss a\sss face of\dss some simplex $\sigma$
entering\trs $\beta$ with non-zero coefficient\sss
and\dss
$\tau\dff \cap\dff X\fff'\off \neq\off \varnothing$\nnsp,\oss
then also\dss
$\sigma\dff \cap\dff X\fff'\off \neq\off \varnothing$\nnsp.\oss
Therefore\qss  (\ref{gamma})\qss implies\sss that\dss\vspace{3pt}
\[
\quad
\gamma\fff'
\off =\off 
r\fff'
\qff +\qff
(\trf \partial\trf \beta\trf)\trf|\qff X\fff'
\pff.
\]

\vspace{-9pt}
Let\dss $F$\sss be\sss the boundary\sss of\dss $X\fff'$\nnsp.\oss

The difference\vspace{3pt}
\[
\quad
d
\off =\off
\partial\trf \beta\fff'
\qff -\qff 
(\trf \partial\trf \beta\trf)\trf|\qff X\fff'
\]

\vspace{-9pt}
is\dss a\sss sum of\trs faces $\tau$ of\dss
simplices $\sigma$
such\dss that\sss
$\tau\dff \cap\trf X\fff'
\off =\off 
\varnothing$\sss 
and\sss 
$\sigma\dff \cap\trf X\fff'
\off \neq\off
\varnothing$\nnsp.\oss
This\sss implies\sss that\dss
$\sigma\dff \cap\trf F
\off \neq\off
\varnothing$\dss
and\dss hence\dss
$\sigma\dff \cap\trf C_{\trf Z}
\off \neq\off
\varnothing$\dss
for some\sss $Z\qff \in\qff \mathcal{Z}$\nnsp.\oss
Together\sss with\sss
$\tau\dff \cap\trf X\fff'
\off =\off 
\varnothing$\sss 
this implies\sss that\sss
$\tau$\sss is\dss a simplex\sss in\dss 
$V_{\trf Z}
\qff \subset\qff
C_{\trf Z}$\nsp.\oss
Since\sss $\beta$\sss is\dss
a compactly\sss finite chain,\oss
it\dss follows\sss that\sss $d\dff \cap\dff C_{\trf Z}$\sss is\dss
a\sss finite chain\sss for every\sss $Z$\sss and\sss 
$d$\sss is\dss a compactly\sss finite chain\sss in\sss $V$\nnsp.\oss
Clearly,\pss\vspace{3pt}
\[
\quad
\gamma\fff'\qff +\qff d
\off =\off
r\fff'\qff +\qff \partial\trf \beta\fff'
\pff.
\]

\vspace{-9pt} 
The construction of\dss $\ry\trf(\dff c\trf)$\sss
shows\sss that\dss there exists a\sss chain $s\fff'$ in\sss $V$\sss
such\dss that\sss 
$r\fff'\qff +\qff s\fff'$\sss is\dss a cycle subdividing\sss $c$\sss
and\sss
$s\fff'\trf|\qff V_{\trf Z}$\sss 
is\dss a\sss finite chain\sss for every\sss 
$Z\qff \in\qff \mathcal{Z}$\nnsp.\oss 
It\dss follows\sss that\sss
$r\fff'\qff +\qff s\fff'$\sss is\dss a cycle 
representing\sss $h$\sss and\sss $s\fff'$\sss
is\dss a compactly\sss finite chain.\oss
Since\sss $d$\sss is\dss also compactly\sss finite,\oss
the chain\sss $s\off =\off d\qff +\qff s\fff'$\sss 
is\dss compactly\sss finite.\oss
Clearly,\vspace{3pt}
\[
\quad
\gamma\fff'\qff +\qff s
\off =\off
\gamma\fff'
\qff +\qff d\qff +\qff s\fff'
\off =\off
r\fff'\qff +\qff s\fff'\qff +\qff \partial\trf \beta\fff'
\pff
\]

\vspace{-9pt}
and\dss hence\sss
$\gamma\fff'\qff +\qff s$\sss 
is\dss a cycle representing\sss $h$\nnsp.\oss
The\sss lemma\sss follows.\oss  \eproof

\mypar{Theorem.}{removing-lone}
\emph{Suppose\sss that\dss $\mathcal{Z}$\sss is\dss countable
and\dss for\sss every\qss
$Z\qff \in\qff \mathcal{Z}$ a set\qss $Z_{\trf +}$\dss is\dss given,\oss
such\dss that\sss $C_{\trf Z}$\sss is\dss weakly $l_{\dff 1}$\dnsp-acyclic\sss in\sss $Z_{\trf +}$\nsp.\oss
If\pss the family\sss of\trs sets\dss $Z_{\trf +}$
is\dss compactly\trs finite,\oss
then}\vspace{3pt}
\[
\quad
\norm{\ry\trf(\dff h\trf)}
\off \geq\off
\norm{h}
\pff
\]

\vspace{-9pt}
\emph{for every\dss homology\sss class\qss 
$h\qff \in\qff H_{\dff n}^\cf\dff(\trf X\trf)$\nnsp.\oss}

\proof
If\qss
$\norm{\ry\trf(\dff h\trf)}\off =\off \infty$\nnsp,\oss
there\dss is\dss nothing\sss to prove.\oss
Suppose\sss that\sss
$\norm{\ry\trf(\dff h\trf)}\qff <\qff \infty$\nnsp.\oss
Let\dss us\sss fix an arbitrary $\varepsilon\qff >\qff 0$\nnsp.\oss
Then\dss there exists a compactly\sss finite cycle\sss 
$\gamma$\sss in\sss $X\ssm Y$\sss
such\dss that\sss\vspace{3pt}
\[
\quad
\norm{\gamma}
\off <\off
\norm{\ry\trf(\dff h\trf)}
\qff +\qff
\varepsilon
\]

\vspace{-9pt}
and\sss $\gamma$\sss represents\sss the homology\sss class\sss $\ry\trf(\dff h\trf)$\nnsp.\oss
Since\sss $\mathcal{Z}$\sss is\dss countable,\oss\vspace{3pt}\vspace{-0.125pt}
\[
\quad
\varepsilon
\off =\off
\sum\nolimits_{\qff Z\qff \in\qff \mathcal{Z}}\qff
\varepsilon_{\dff Z}
\]

\vspace{-9pt}
for some numbers\sss
$\varepsilon_{\dff Z}\qff >\qff 0$\nnsp.\oss
For every\dss $Z\qff \in\qff \mathcal{Z}$\qss let\dss
$\delta_{\trf Z}$\sss be some number\dss $>\qff 0$\dss 
such\dss that\dss the
conclusion of\qss Lemma\qss \ref{double-boundary}\qss
holds\sss for\dss\vspace{3pt}
\[
\quad
C_{\trf Z}\dff,\quad Z_{\trf +}\dff,\quad \varepsilon_{\dff Z}\dff,\quad
\mbox{and}\quad
\delta_{\trf Z}
\]

\vspace{-9pt}
in\sss the roles\sss of\pss
$C\fff,\off\dff D\fff,\off\dff \varepsilon$\dss 
and\dss $\delta$\dss respectively.\oss\vspace{1.125pt}

By\qss Lemma\qss \ref{chain-modification}\qss for\sss every\sss $Z$\sss
there exists an open\sss neighborhood\dss 
$U_{\dff Z}$\sss of\trs $Z$\dss in\sss $C_{\trf Z}$\sss
such\dss that\vspace{3pt}
\begin{equation}
\label{z-delta}
\quad
\norm{(\trf \partial\trf (\trf \gamma\dff \cap\dff U_{\trf Z}\trf)}
\off <\off\qff
\delta_{\trf Z}
\pff.
\end{equation}

\vspace{-9pt}
Let\sss $U$\sss be\sss the union of\trs the neighborhoods\sss $U_{\trf Z}$\sss
and\dss let\dss
$X\fff'\off =\off X\ssm U$\nnsp.\oss
By\qss Lemma\qss \ref{homology}\qss 
there exists\sss a compactly\trs finite chain\sss $s$\sss in\dss $V$\sss
such\dss that\dss 
$\gamma\trf |\qff X\fff'\pff +\pff s$\dss 
is\dss a cycle representing\sss $h$\nnsp.\oss
Since\sss $V_{\dff Z}\dff,\pff U_{\trf Z}\qff \subset\qff C_{\trf Z}$\sss
and\dss the sets\sss $C_{\trf Z}$\sss are pair-wise disjoint,\oss\vspace{4.5pt}
\[
\quad
\gamma\dff \cap\dff U
\off =\off
\sum\nolimits_{\qff Z\qff \in\qff \mathcal{Z}}\qff
\gamma\dff \cap\dff U_{\trf Z}
\quad\
\mbox{and}\quad\
s\dff \cap\dff V
\off =\off
\sum\nolimits_{\qff Z\qff \in\qff \mathcal{Z}}\qff
s\dff \cap\dff V_{\trf Z}
\pff.
\]

\vspace{-7.5pt}
Since\dss 
$\gamma\trf |\qff X\fff'\pff +\pff s$\dss 
and\dss
$\gamma
\off =\off
\gamma\trf |\qff X\fff'
\pff +\pff
\gamma\dff \cap\dff U$\sss
are cycles,\oss\vspace{3pt}
\[
\quad
\partial\trf (\trf \gamma\dff \cap\dff U\trf)
\off =\off
-\qff
\partial\trf (\trf \gamma\trf |\qff X\fff'\trf)
\off =\off
\partial\trf (\trf s\dff \cap\dff V\trf)
\]

\vspace{-9pt}
and\dss hence\qss
$\partial\trf (\trf \gamma\dff \cap\dff U_{\trf Z}\trf)
\off =\off
\partial\trf (\trf s\dff \cap\dff V_{\trf Z}\trf)$\qss 
for\sss every\sss $Z$\nnsp.\oss
In\sss view of\trs the choice of\dss $\delta_{\trf Z}$\nsp,\oss
the inequality\qss (\ref{z-delta})\qss together\sss with\trs
Lemma\qss \ref{double-boundary}\qss imply\dss that\dss
for\sss every\sss $Z$\sss
there exists a\sss finite chain\sss $\psi_{\trf Z}$\sss in\sss $Z_{\trf +}$\sss
such\dss that\sss
$\psi_{\trf Z}\qff -\qff (\trf s\dff \cap\dff V_{\trf Z}\trf)$\sss 
is\dss the boundary\sss of\dss a\sss finite chain\sss $b_{\trf Z}$\dss
in\dss $Z_{\trf +}$\sss and\dss 
$\norm{\psi_{\trf Z}}\qff <\pff \varepsilon_{\trf Z}$\nsp.\oss
Let\vspace{4.5pt}
\[
\quad
\psi
\off =\off
\sum\nolimits_{\qff Z\qff \in\qff \mathcal{Z}}\qff
\psi_{\trf Z}
\quad\
\mbox{and}\quad\
b
\off =\off
\sum\nolimits_{\qff Z\qff \in\qff \mathcal{Z}}\qff
b_{\trf Z}
\pff.
\]

\vspace{-7.5pt}
Then\sss
$\psi\qff -\qff s
\off =\off
\partial\dff b$\sss
and\dss hence\vspace{3pt}
\[
\quad
\gamma\trf |\qff X\fff'\pff +\pff \psi
\off =\off
\gamma\trf |\qff X\fff'\pff +\pff s
\pff +\pff
(\trf \psi\qff -\qff s\trf)
\off =\off
\gamma\trf |\qff X\fff'\pff +\pff s
\pff +\pff
\partial\dff b
\pff.
\]

\vspace{-9pt}
Since\sss the family\sss of\trs the sets\sss $Z_{\trf +}$\sss is\dss
compactly\sss finite,\pss $\psi$\sss and\dss $b$\sss
are compactly\sss finite chains.\oss
It\dss follows\sss that\sss
$\gamma\trf |\qff X\fff'\pff +\pff \psi$\sss
is\dss a compactly\sss finite cycle representing\sss $h$\nnsp.\oss
At\dss the same\sss time\vspace{4.5pt}
\[
\quad
\norm{\psi}
\off \leq\dff\off
\sum\nolimits_{\qff Z\qff \in\qff \mathcal{Z}}\qff
\norm{\psi_{\trf Z}}
\off \leq\dff\off
\sum\nolimits_{\qff Z\qff \in\qff \mathcal{Z}}\qff
\varepsilon_{\dff Z}
\off =\off
\varepsilon
\]

\vspace{-12pt}
and\dss hence\vspace{0pt}
\[
\quad
\norm{\gamma\trf |\qff X\fff'
\qff +\qff
\psi}
\off \leq\off
\norm{\gamma\trf |\qff X\fff'}
\qff +\qff
\norm{\psi}
\]

\vspace{-36pt}
\[
\quad
\phantom{\norm{\gamma\trf |\qff X\fff'
\qff +\qff
\psi}
\off }
\leq\off
\norm{\gamma\trf |\qff X\fff'}
\qff +\qff
\varepsilon
\off \leq\off
\norm{\gamma}
\qff +\qff
\varepsilon
\]

\vspace{-36pt}
\[
\quad
\phantom{\norm{\gamma\trf |\qff X\fff'
\qff +\qff
\psi}
\off <\off
\norm{\gamma\trf |\qff X\fff'}
\qff +\qff
\varepsilon
\off }
<\off
\norm{\ry\trf(\dff h\trf)}
\qff +\qff
\varepsilon
\qff +\qff
\varepsilon
\off =\off
\norm{\ry\trf(\dff h\trf)}
\qff +\qff
2\dff \varepsilon
\pff.
\]

\vspace{-7.5pt}
Since\sss $\varepsilon\qff >\qff 0$\sss is\dss arbitrary,\oss
it\dss follows\dss that\sss
$\norm{h}
\off \leq\off
\norm{\ry\trf(\dff h\trf)}$\nnsp.\oss  \eproof

\myuppar{Gromov's\trs Cutting-of\qss theorem.}
In\dss this\sss theorem\qss 
(see\qss \cite{gro},\oss Theorem\qss (2)\qss in\dss Section\qss 4.2)\qss
Gromov\dss assumes\sss that\sss $X$\sss is\dss a manifold,\pss
$Y$\sss is\dss the union of\dss a sequence of\dss 
disjoint\sss compact\dss submanifolds\sss $Y_{\fff i}$\nsp,\qss $i\qff \in\pff \nnn$\nnsp,\oss
possibly\sss with\sss boundary,\oss
every\sss $Y_{\fff i}$\sss is\dss amenable\sss in\sss $X$\sss in\sss
the sense of\qss Gromov,\oss
and\dss the family\sss 
$\{\trf Y_{\fff i}\trf\}_{\qff i\dff \in\dff \nnn}$\sss is\dss 
almost\sss compactly\sss 
amenable in\sss the sense of\qss Gromov\trs in\sss $X$\dss
(see\dss Section\qss \ref{extensions}\qss for\dss the definitions).\oss
Since a compact\dss manifold can\sss have only\sss finitely\sss components,\oss
under\sss these assumptions\sss the family\sss of\dss components of\dss $Y$\sss
is\dss compactly\sss amenable and\dss hence\dss is\dss compactly $l_{\dff 1}$\dnsp-acyclic.\oss
Since components\sss $Z$\sss of\dss $Y$\sss are submanifolds,\oss
standard\dss results imply\dss the existence of\dss compact\dss neighborhoods\sss
$C_{\trf Z}$\sss with\sss the properties required above.\oss
Therefore\trs Theorem\qss \ref{removing-lone}\qss applies\sss under\trs
Gromov's\dss assumptions.\oss
Its conclusion\dss is\dss the same as\trs Gromov's\qss inequality\sss
$\norm{h\fff'}\qff \geq\qff \norm{h}$\nnsp.\oss

\myappend{Double\qss complexes}{double-complexes}

\myuppar{Double complexes.}
In order\sss to deal\dss with almost\sss $l_{\dff 1}$\dnsp-acyclic
coverings,\oss
we need\sss a complement\dss to\sss the\sss theorem about\sss
double complexes.\oss
Let\sss $K_{\dff \bullet\fff,\dff \bullet}$\sss
be a double complex with differentials\sss
$d\dff \colon\dff
K_{\dff p\fff,\dff q}
\qff \ttoo\qff
K_{\dff p\fff,\dff q\dff -\dff 1}$\sss
and\dss
$\delta\dff \colon\dff
K_{\dff p\fff,\dff q}
\qff \ttoo\qff
K_{\dff p\dff -\dff 1\fff,\dff q}$\nnsp.\oss
Let\sss $T_{\dff \bullet}$\sss be\sss the\sss total\sss complex of\dss
$K_{\trf \bullet\fff,\dff \bullet}$\dnsp,\oss
and\dss $E_{\dff p}$\sss be the cokernel of\qss
$d\dff \colon\dff
K_{\trf p\fff,\dff 1}\qff \ttoo\qff K_{\trf p\fff,\dff 0}$\dnsp,\oss
i.e.\vspace{3pt}\vspace{-1.5pt}
\[
\quad
E_{\dff p}
\off =\off\dff
K_{\trf p\fff,\dff 0}\dff\bigl/d\dff\left(\qff K_{\trf p\fff,\dff 1} \qff\right)
\pff.
\]

\vspace{-9pt}\vspace{-1.5pt}
Recall\dss that\sss $\delta$\sss induces homomorphisms\qss
$\delta_{\trf E}\qff \colon\qff
E_{\dff p}
\qff \ttoo\qff
E_{\dff p\dff -\dff 1}$\qss
turning\dss $E_{\dff \bullet}$\dss into\sss a\sss complex.\oss
Clear\-ly,\pss $E_{\dff \bullet}$\sss is\dss 
a\sss quotient\sss of\dss $T_{\dff \bullet}$\nsp,\oss
and\dss $K_{\trf 0\fff,\dff \bullet}$\sss with\dss the differential\sss $d$\sss
is\dss a subcomplex of\dss $T_{\dff \bullet}$\nsp.\oss

\myapar{Lemma.}{double-complex}
\emph{If\trs the complexes\dss
$(\qff K_{\dff p\fff,\dff \bullet}\fff,\qff d\qff)$\sss
with\dss $p\qff >\qff 0$\sss are exact,\oss
then\dss the\sss kernel\sss of\qss the map\dss
$H_{\dff \bullet}\dff(\trf T_{\dff \bullet}\trf)
\qff \ttoo\qff
H_{\dff \bullet}\dff(\trf E_{\dff \bullet}\trf)$\dss
is\dss contained\sss in\dss the image of\qss the map\dss
$H_{\dff \bullet}\dff(\trf K_{\trf 0\fff,\dff \bullet}\fff,\qff d\trf)
\qff \ttoo\qff
H_{\dff \bullet}\dff(\trf T_{\dff \bullet}\trf)$\nnsp.\oss}

\proof
The\sss proof\trs is\dss completely\sss similar\sss to\sss
the proof\dss of\dss injectivity\sss in\sss the proof\dss
of\qss Theorem\qss A.2\qss in\qss \cite{i3}.\oss
We\sss leave details\sss to\sss the reader\halfff.\oss  \eproof

\begin{flushright}

December 14\fff,\oss 2020
 
https\halfff:/\!/\!nikolaivivanov.com

E-mail\halfff:\oss nikolai.v.ivanov{\fff}@{\dff}icloud.com

\end{flushright}

\end{document}